\documentclass[titlepage,final,11pt]{article}
\usepackage{theorem}
\usepackage{enumerate}
\usepackage{array}
\usepackage[reqno]{amsmath}
\usepackage{amssymb}
\usepackage{mathtools}
\usepackage{latexsym}
\usepackage{makeidx}
\usepackage{fancybox}
\input epsf.sty
\usepackage[dvips]{graphicx,color}
\usepackage{showlabels}
\usepackage{subcaption}
\usepackage[dvipsnames]{xcolor}
\usepackage[normalem]{ulem}

\theoremstyle{change}
\newtheorem{proclaim}{PROCLAIM}[section]
\newtheorem{theorem}[proclaim]{Theorem}
\newtheorem{definition}[proclaim]{Definition}

\newtheorem{proposition}[proclaim]{Proposition}
\newtheorem{corollary}[proclaim]{Corollary}
\newtheorem{example}[proclaim]{Example}

\numberwithin{equation}{section}

\outer\def\proclaim #1. #2\par{\medbreak \noindent{\bf#1.\enspace}{\sl#2}\par
  \ifdim\lastskip<\medskipamount
  \removelastskip\penalty55\medskip\fi}
\def\state #1. { \noindent{\bf#1.\enspace}}
\def\algo #1. { \noindent{\bf#1.\enspace}}

\DeclareMathOperator*{\argmin}{argmin}

\DeclareMathOperator{\dist}{dist}

\DeclareMathOperator{\dom}{dom}
\DeclareMathOperator{\gph}{gph}
\DeclareMathOperator{\epi}{epi}

\newcommand{\comp}{\,{\raise 1pt \hbox{$\scriptstyle\circ$}}\,}

\newcommand{\reals}{\mathbb{R}}
\newcommand{\Reals}{\overline{\mathbb{R}}}
\newcommand{\natnums}{{{\rm l} \kern -.13em {\rm N} }}
\newcommand{\nats}{\mathbb{N}}
\newcommand{\snats}{{I\kern -.29em N}}
\newcommand{\rats}{{Q\kern -.64em \raise 1pt \hbox{$\scriptstyle |$}\;\,}}
\newcommand{\srats}
	{{Q\kern -.56em \raise 1.2pt \hbox{$\scriptscriptstyle /$}\,}}
\newcommand{\ints}{Z\kern -.46em Z}
\newcommand{\ball}{\mathbb{B}}
\newcommand{\pluss}{\hskip1pt \raise1pt\vbox{\hrule width6pt \vskip1pt \hrule
                    width6pt} \kern-4pt{\lower1pt\hbox{\vrule height6pt
		    \kern1pt\vrule height6pt}}\hskip5pt}
\newcommand{\eop}
	{\hfill{$\vcenter{\hrule height1pt \hbox{\vrule width1pt height5pt
   	 \kern5pt \vrule width1pt} \hrule height1pt}$} \medskip}
\newcommand{\half}
	{{\raisebox{1pt}{$\frac{1}{2}$}}}

\newcommand{\lminus}{{\scriptscriptstyle -}}
\newcommand{\lplus}{{\scriptscriptstyle +}}

\newcommand{\setd}{{ d \kern -.15em l}}
\newcommand{\hatsetd}{ d \hat{\kern -.15em l }}

\renewcommand{\epsilon}{\varepsilon}
\renewcommand{\phi}{\varphi}

\newcommand{\tto}{\;{\lower 1pt \hbox{$\rightarrow$}}\kern -12pt
           \hbox{\raise 2.5pt \hbox{$\rightarrow$}}\;}
\newcommand{\overto}[1]{\,{\raise 0pt\hbox{$\rightarrow$}}\kern -9pt
     \hbox{\lower 3pt \hbox{$\scriptscriptstyle#1$}}\hskip6pt}
\newcommand{\underto}[1]{\,{\lower 1pt\hbox{$\rightarrow$}}\kern -9pt
     \hbox{\raise 4pt \hbox{$\,\scriptscriptstyle#1$}}\hskip7pt}
\newcommand{\bigoverto}[1]{{\raise 0pt\hbox{$\,\longrightarrow$}}\kern -16pt
     \hbox{\lower 3pt \hbox{$\scriptscriptstyle#1$}}\hskip4pt}
\newcommand{\bigunderto}[1]{\,{\lower 1pt\hbox{$\longrightarrow$}}\kern -16pt
     \hbox{\raise 4pt \hbox{$\,\scriptscriptstyle#1$}}\hskip6pt}
\newcommand{\bigbigto}[2]{\,{\raise 0pt\hbox{$\,\longrightarrow$}}\kern -16pt
     \hbox{\lower 3pt \hbox{$\scriptscriptstyle#2$}}\kern -10pt
     \hbox{\raise 4pt \hbox{$\,\scriptscriptstyle#1$}}\hskip7pt}
\newcommand{\downto}{{\raise 1pt \hbox{$\scriptscriptstyle \,\searrow\,$}}}
\newcommand{\upto}{{\raise 1pt \hbox{$\scriptscriptstyle \,\nearrow\,$}}}

\newcommand{\notimply}
	{\quad\hbox{$\Longrightarrow \kern -14pt {/}$}\hskip6pt\quad}

\newcommand{\lto}{\,{\lower 1pt\hbox{$\rightarrow$}}\kern -10pt
     \hbox{\raise 4pt \hbox{$\, \scriptstyle l$}}\hskip7pt}
\newcommand{\eto}{\,{\lower 1pt\hbox{$\rightarrow$}}\kern -10pt
     \hbox{\raise 4pt \hbox{$\, \scriptstyle e$}}\hskip7pt}
\newcommand{\hto}{\,{\lower 1pt\hbox{$\rightarrow$}}\kern -11pt
     \hbox{\raise 4pt \hbox{$\, \scriptstyle h$}}\hskip7pt}
\newcommand{\pto}{\,{\lower 1pt\hbox{$\rightarrow$}}\kern -11pt
     \hbox{\raise 4.5pt \hbox{$\, \scriptstyle p$}}\hskip7pt}
\newcommand{\cto}{\,{\lower 1pt\hbox{$\rightarrow$}}\kern -11pt
     \hbox{\raise 4pt \hbox{$\, \scriptstyle c$}}\hskip7pt}
\newcommand{\gto}{\,{\lower 1pt\hbox{$\rightarrow$}}\kern -11pt
     \hbox{\raise 4.5pt \hbox{$\, \scriptstyle g$}}\hskip7pt}
\newcommand{\sto}{\,{\lower 1pt\hbox{$\rightarrow$}}\kern -10pt
     \hbox{\raise 4pt \hbox{$\, \scriptstyle s$}}\hskip7pt}
\newcommand{\awto}{\,{\lower 1pt\hbox{$\rightarrow$}}\kern -15pt
     \hbox{\raise 4pt \hbox{$\, \scriptstyle aw$}}\hskip7pt}
\def\Nto{\,{\raise 1pt\hbox{$\rightarrow$}}\kern -13pt
     \hbox{\lower 3pt \hbox{$\, \scriptstyle N$}}\hskip7pt}
\def\Cto{\,{\raise 1pt\hbox{$\rightarrow$}}\kern -14pt
     \hbox{\lower 3pt \hbox{$\, \scriptstyle C$}}\hskip7pt}
\def\fto{\,{\raise 1pt\hbox{$\rightarrow$}}\kern -14pt
     \hbox{\lower 3pt \hbox{$\, \scriptstyle f$}}\hskip7pt}

\newcommand{\low}[1]{{\lower1pt \hbox{$\scriptstyle #1$}}}
\newcommand{\loww}[1]{{\lower2pt \hbox{$\scriptstyle #1$}}}
\newcommand{\high}[1]{{\raise1pt \hbox{$\scriptstyle #1$}}}

\newcommand{\cV}{{\cal V}}

\newcommand{\nsum}{\mathop{\sum}\nolimits}

\newcommand{\nOutLim}{\mathop{\rm Lim\hspace{-0.01cm}Out}\nolimits}

\newcommand{\nlim}{\mathop{\rm lim}\nolimits}
\newcommand{\nliminf}{\mathop{\rm liminf}\nolimits}
\renewcommand{\liminf}{\mathop{\rm liminf}}
\renewcommand{\limsup}{\mathop{\rm limsup}}
\newcommand{\nlimsup}{\mathop{\rm limsup}\nolimits}
\newcommand{\ninf}{\mathop{\rm inf}\nolimits}
\newcommand{\nsup}{\mathop{\rm sup}\nolimits}

\newcommand{\nmin}{\mathop{\rm min}\nolimits}
\newcommand{\nnmin}{\mathop{\rm minimize}}

\newcommand{\nargmin}{\mathop{\rm argmin}\nolimits}

\newcommand{\prob}{{\mathop{\rm prob}\nolimits}}

\newcommand{\bfxi}{\mbox{\boldmath $\xi$}}

\newcommand{\bfy}{\mbox{\boldmath $y$}}

\newcommand{\bfz}{\mbox{\boldmath $z$}}

\newcommand{\bfx}{\mbox{\boldmath $x$}}

\newcommand{\lwdy}[2]{\mathrel{\mathop
        {\raisebox{0.1ex}{\null$#1$}}{\hbox{\kern -1.0em
	{\raisebox{-0.8ex}{$\scriptstyle{\;\to #2}$}}}}}}
\newcommand{\lwwdy}[2]{\mathrel{\mathop
        {\raisebox{0.2ex}{\null$#1$}}{\hbox{\kern -1.0em
	{\raisebox{-1.1ex}{$\scriptstyle{\;\to #2}$}}}}}}
\newcommand{\slwdy}[2]{\scriptsize{{\mathrel{\mathop
        {\raisebox{0.1ex}{\null$#1$}}{\hbox{\kern -1.0em
	{\raisebox{-0.8ex}{$\scriptstyle{\;\to #2}$}}}}}}}}
\newcommand{\slwwdy}[2]{\scriptsize{{\mathrel{\mathop
        {\raisebox{0.2ex}{\null$#1$}}{\hbox{\kern -1.0em
	{\raisebox{-1.1ex}{$\scriptstyle{\;\to #2}$}}}}}}}}

\definecolor{lightgray}{gray}{0.75}
\definecolor{myred}{rgb}{0.55,0,0}
\definecolor{myblue}{rgb}{0,0,0.5} 
\definecolor{mygreen}{rgb}{0,0.5,0} 
\definecolor{purple}{rgb}{0.5,0,0.5} 
\definecolor{turq}{rgb}{0,0.805,0.816} 
\definecolor{maroon}{rgb}{0.51,0,0}
\definecolor{MAROON}{rgb}{0.51,0,0}
\definecolor{redor}{rgb}{0.78,0.078,0.078}
\definecolor{dgreen}{rgb}{0,0.3,0}

\newcommand{\Ex}{\mathbb{E}}

\newcommand{\bcdot}{\,{\raise .2ex \hbox{$\centerdot$}}\,}

\usepackage{amsmath}
\usepackage[toc,page]{appendix}
\mathchardef\mhyphen="2D 

\usepackage[normalem]{ulem}

\numberwithin{subcase}{case}

\begin{document}


\begin{center}
\begin{large}
{\bf Rockafellian Relaxation and Stochastic Optimization under Perturbations}
\smallskip
\end{large}
\vglue 0.3truecm
\begin{tabular}{c}
  \begin{large} {\sl Johannes O. Royset, ~~ Louis L. Chen, ~~ Eric Eckstrand 
                                  } \end{large} \\
  Operations Research Department\\
  Naval Postgraduate School, Monterey, California\\
\end{tabular}

\vskip 0.1truecm

{\bf Date}:\quad \ \today

\end{center}

\vskip 0.3truecm

\noindent {\bf Abstract}.
In practice, optimization models are often prone to unavoidable inaccuracies due to dubious assumptions and corrupted data. Traditionally, this placed special emphasis on risk-based and robust formulations, and their focus on ``conservative" decisions. We develop, in contrast, an ``optimistic" framework based on Rockafellian relaxations in which optimization is conducted not only over the original decision space but also jointly with a choice of model perturbation. The framework enables us to address challenging problems with ambiguous probability distributions from the areas of two-stage stochastic optimization without relatively complete recourse, probability functions lacking continuity properties, expectation constraints, and outlier analysis. We are also able to circumvent the fundamental difficulty in stochastic optimization that convergence of distributions fails to guarantee convergence of expectations. The framework centers on the novel concepts of exact and limit-exact Rockafellians, with interpretations of ``negative'' regularization emerging in certain settings. We illustrate the role of Phi-divergence, examine rates of convergence under changing distributions, and explore extensions to first-order optimality conditions. The main development is free of assumptions about convexity, smoothness, and even continuity of objective functions. Numerical results in the setting of computer vision and text analytics with label noise illustrate the framework.

\vskip 0.2truecm


\baselineskip=15pt

\section{Introduction}

Uncertain parameter values in stochastic optimization are modeled using probability distributions, but these can themselves be unsettled due to a lack of data or incorrect information. Even the slightest error in the prescribed distribution could make a truly optimal solution look poor. This vulnerability also underscores the threat of adversarially-caused corruption, which could be strategically crafted in advance to mislead naive decision-makers away from good solutions. In this work, we propose a framework that supports decision-making robust to modeling errors, adversarial perturbations, and other deviations from an actual problem. Decision-makers are oftentimes unaware of the ``true" form to a problem so that they cannot detect the presence of deviations, much less determine how such deviations can be corrected. The best that could be hoped for is that their decisions made based on solving slightly deviated problems are also only slightly deviated from those of the unknown actual problem. We aim to arm the decision-maker with a framework that exhibits this kind of continuity properties in the setting of expectation functions.

Risk-based and distributionally robust approaches to minimizing expectation functions furnish performance guarantees across a collection of possible distributions, among which we hope to find the actual one. Thus, it replaces expectation functions by higher, more conservative approximations. It is easily seen via Fatou's lemma that in some applications the move to ``high'' approximations is misguided: A weakly converging sequence of nonnegative random variables can have expectations much above that of the limiting random variable. Consequently, even a nearly correct distribution may already produce a value too high. In this paper, we address such situations by constructing ``low'' approximations via {\em Rockafellian relaxations} and thereby consider best-case instead of worst-case distributions. The focus on relaxation over restriction of a problem allows us to keep in play decisions that {\em might} be good under the actual probability distribution.

Undoubtedly, problem relaxations lead to decisions that could be less conservative but may also uncover previously hidden possibilities. This is especially the case for tightly constrained problems, for problems with induced constraints such as in the absence of relatively complete recourse, for problems with discontinuous functions expressing quantiles and probabilities, and for statistical problems with outliers and other abnormalities in the training data. Through examples, we show that even in the case of finite probability spaces an arbitrarily accurate estimate of the actual distribution can cause large errors in optimal solutions. We develop a framework for addressing such situations and examine convergence as distributions become more accurate. In some cases, the framework has the interpretation of inducing a negative regularization to a problem.

The concept of Rockafellian relaxation can be traced back to \cite{Rockafellar.63} and the perspective of viewing a minimization problem as a member of a family of perturbed problems, with this family being specified by a function of both decision variables and perturbation parameters. We refer to such functions as {\em Rockafellians}, while the name ``bifunctions'' appears in the pioneering development of convex duality  \cite[Chapter 29]{Rockafellar.70} and also in extensions to infinite-dimensional convex analysis \cite{Rockafellar.74} and nonconvex problems \cite{Rockafellar.85}. The name ``Rockafellian'' emerges in \cite[Chapter 5]{primer} and \cite{Royset.21}, with ``perturbation function'' \cite{zalinescu2002convex} and ``bivariate function'' \cite{bauschke2011convex} also being found in the literature.

The perturbation scheme defined by a Rockafellian produces a relaxation of the actual problem involving minimization over both decision variables and perturbation parameters. If designed and tuned appropriately, however, the relaxation may turn out to be equivalent to the actual problem in some sense. The tasks of tuning relaxations are traditionally referred to as dual problems, with classical linear programming and augmented Lagrangian duality emerging from specific choices of Rockafellians; see \cite[Chapters 5 and 6]{primer} and \cite[Chapter 11]{VaAn}. We focus on perturbations that stem from ambiguity about probability distributions in stochastic optimization. This brings us beyond classical perturbation schemes, with tuning becoming less critical. In fact, our development is related to penalty methods and especially exact penalization in constrained optimization (see, e.g., \cite{burke1991exact}). This is reflected in our choice of terminology: {\em exact} and {\em limit-exact} Rockafellians emerge here for the first time.

The systematic study of stochastic optimization problems under changes to probability distributions extends at least back to \cite{RobinsonWets.87,Kall.87,Dupacova.90,RomischSchultz.91}; see \cite{RomischWets.07} and \cite{DentchevaRomisch.13} for more recent efforts in the context of convex problems and stochastic dominance, respectively. There is an extensive literature on how to address parameter uncertainty in a {\em conservative} manner; see \cite{polak1997optimization,stein2003bi,royset2012rate} for results from semi-infinite programming and \cite{ben2002robust,ben2009robust,bertsimas2011theory} for robust convex optimization. Under distributional ambiguity, the focus tends to be on achieving a solution that is ``good'' across a set of candidate distributions, for example centered at a presently available empirical distribution \cite{ben2013robust,WiesemannKuhnSim.14,RoysetWets.16b,mohajerin2018data,bertsimas2018data,DuchiGlynnNamkoong.21}. Related risk-based approaches have similar effect, at least for monotone and positively homogeneous regular measures of risk \cite{FollmerSchied.04,RockafellarUryasev.13,RockafellarRoyset.15}. Adversarial statistical learning addresses ambiguity about training data through conservative perturbations of the support of the corresponding (empirical) distribution \cite{Madry.18}.
Robust regularization \cite{LewisPang.10} achieves conservativeness through perturbation of a decision vector, with applications to manufacturing \cite{MenFreundNguyenSaaseoanePeraire.14} and statistical learning \cite{wu2020adversarial,TsaiHsuYuChen.21,NortonRoyset.22} in the form of diametrical risk minimization. When approximations stem from empirical distributions, these approaches are well-motivated by the downward bias in sample average approximations \cite{NortonRoyset.22}. We step away from conservative approaches, which correspond to being ``pessimistic'' about values of ambiguous parameters and distributions, and adopt an ``optimistic'' approach based on problem relaxation.

The optimistic perspective has received less attention, with much of it centered on statistical learning. In \cite{norton2017optimistic}, both pessimistic and optimistic disturbances are simultaneously considered to form optimization problems that are nonconvex, but nonetheless have structure that remains computationally tractable. This unpublished manuscript notes that nonconvex regularization and other practices for dealing with outliers have connections to such formulations; see also the recent preprint \cite{jiangdfo} for further connections with robust statistics. In \cite{nguyen2019optimistic}, the task of likelihood approximation in a Bayesian setting is framed as problems also involving optimization over probability measures constrained by KL-divergences, moment constraints, or Wasserstein distances. Optimistic formulations with KL-divergence penalties can sometimes outperform solutions from sample average approximations \cite{gotoh2021data}. Additional optimistic formulations emerge in reinforcement learning \cite{song2020optimistic,agarwal2020optimistic}, in decision-making involving covariate information \cite{cao2021contextual}, in outlier analysis \cite{AravkinDavis.19,Zhengetal.21,Narasimhanetal.23}, and in chance-constrained models \cite{HanasusantoRoitchKuhnWiesemann.17}.

There are connections between pessimistic formulations of a primal problem and optimistic formulations of a corresponding dual problem \cite{beck2009duality}; see also \cite{jeyakumar2010strong,li2011robust,suzuki2013surrogate} for such insight in the context of convex problems. These duality relationships extend to distributionally robust optimization \cite{zhen2021mathematical}.

The present paper stands out from the literature on optimistic problem formulations by adopting exceptionally mild assumptions. We rely mostly on functions being extended real-valued, proper, and lower semicontinuous. In such general settings, optimistic formulations obtained via Rockafellian relaxation prove to be well behaved under changes to probability distributions. We establish a pathway for constructing suitable optimistic formulations. It originates with perturbation theory for optimization problems and facilitates the analysis of many settings beyond those specifically considered here. While extensions are possible, we focus on finite probability spaces to avoid technical distractions; a subsection indicates additional possibilities without this limitation.

We proceed in Section 2 by laying out examples from statistical learning and stochastic optimization. Section 3 develops the general framework of Rockafellian relaxation and the notions of exactness and limiting exactness. Section 4 demonstrates the framework in the context of stochastic optimization. Section 5 turns to optimality conditions stemming from Rockafellian relaxation. Section 6 presents numerical results. An appendix gives additional proofs.

\section{Challenges and Examples} \label{section:: Examples}

Numerous applications involve expectation functions of the form
\[
x\mapsto \Ex_{P}\big[g(\bfxi,x)\big] = \int g(\xi,x) dP(\xi),~ \mbox{ where }~ g:\reals^m\times \reals^n \to \Reals = [-\infty,\infty],
\]
describes a quantity of interest depending on a decision vector $x$ and a random vector $\bfxi$ with probability distribution $P$. (We express random quantities using boldface letters and return to a regular font for their realizations.) As long as $g(\cdot,x)$ is measurable for all $x$, the expectation function is well defined but possible not finite \cite[Section 3.B]{primer} under the usual conventions of extended real-valued arithmetic which we assume throughout\footnote{Most significantly, $0\cdot\infty = 0$ and  $\infty - \infty=\infty$; see \cite[Section 1.D]{primer}.}.

In practice, the actual probability distribution $P$ of $\bfxi$ is rarely known and we are only in possession of an approximating probability distribution obtained through modeling and data analysis, a process that might be affected by adversarial interference. It is therefore critical to examine the effect of replacing $P$ by an alternative on the solutions of the corresponding optimization problems. With this goal in mind, we consider a sequence of probability distributions $\{P^\nu, \nu\in \nats\}$, where $\nats = \{1, 2, 3, \dots\}$, tending to $P$ in some sense as $\nu \to \infty$. The sequence of distributions may not be the outcome of a specific estimation procedure and could be rather conceptual as part of a sensitivity analysis or outlier analysis; the study of influence functions in statistical learning \cite{koh2017understanding} furnishes examples. In particular, $P^\nu$ may {\em not} be an empirical distribution using $\nu$ sample points, although this is a possibility within our framework. Regardless of the setting, one would hope that $P^\nu$ tending to $P$ ensures that optimization problems involving $P^\nu$ indeed produce arbitrarily accurate solutions of the actual problem defined by $P$. It is well known that caution is advised before jumping to such a conclusion.

\begin{example}{\rm (failure of convergence of solutions).}\label{eConvDistr} For $g(\xi,x) = \xi x + \tfrac{1}{2}(1-x)$ and a random variable $\bfxi$ with distribution $P$ assigning  probability $1$ to $\xi = 0$, we obtain that the problem of minimizing $\Ex_{P} [g(\bfxi,x)] = \tfrac{1}{2}(1-x)$ over $x\in [0,1]$ has $x = 1$ as its unique minimizer. Let $P^\nu$ assign probability $1-1/\nu$ to $\xi = 0$ and probability $1/\nu$ to $\xi = \nu$. Then, $P^\nu$ converges weakly to $P$, but minimizing $\Ex_{P^\nu}[g(\bfxi,x)] = \tfrac{1}{2}(1+x)$ over $x \in [0,1]$
has $x = 0$ as its unique minimizer regardless of $\nu\in\nats$.
\end{example}

The example shows that a seemingly accurate $P^\nu$ produces a poor solution relative to the actual problem, with $P^\nu$ leading to an objective function that is too high. This can also occur in the application area of classification under fairness constraints as illustrated next.

\begin{example}{\rm (fair support vector machine).} \label{eFairSVM}
For an affine classifier and a random vector $(\bfx, \bfy, \bfz)$, with $\bfx$ representing features, $\bfy$ giving binary labels, and $\bfz$ specifying sensitive attributes, we consider
\[
\nnmin_{a\in\reals^n, \alpha\in\reals} \,\Ex_{P}\Big[\max\big\{0, 1 - \bfy\big(\langle a, \bfx \rangle + \alpha\big)\big\}\Big] + \|a\|^2_2~ \mbox{ subject to }~ \Ex_{P}\Big[(\bfz - \bar{z})\big(\langle a, \bfx \rangle + \alpha\big)\Big] \leq t,
\]
where $\bar z = \Ex[\bfz]$; see \cite{zafar2017fairness}. The constraint asks for low covariance between $\bfz$ and the classifier. For $n = 1$, $t = 0.25$,  and a distribution $P$ of $(\bfx, \bfy, \bfz)$ assigning probability $1/2$ to $(-1, -1, 0)$ and probability $1/2$ to $(1,1,1)$, we obtain the set of minimizers $\{1/2\} \times [-1/2, 1/2]$. However, if the distribution is slightly altered to assigning probability $1/2$ to $(-1, -1, 0)$, probability $1/2 - 1/\nu$ to $(1,1,1)$, and probability $1/\nu$ to $(\nu, 1, 1)$, then the corresponding minimizers are unique for sufficiently large $\nu$ and they converge to the incorrect $(1/4, -3/4)$. This takes place even though the approximating probability distribution converges weakly to the actual distribution $P$.
\end{example}

Problems with expectation constraints, as in the example and in Neyman-Pearson classification \cite{rigollet2011neyman}, are especially susceptible to changes to probability distributions. As is well known from nonlinear programming, a perturbation of the right-hand side in an inequality constraint can cause large changes in solutions when constraint qualifications fail; see our treatment in Subsection \ref{subsec:constraints}.

Challenges emerge even in the case of a finite probability space. Suppose that $p = (p_1, \dots, p_s)$ is a probability vector, i.e., $p\in \Delta$, with
\[
\Delta = \Big\{ q\in [0, \infty)^s~\Big|~\nsum_{i=1}^s q_i = 1\Big\}.
\]
Then, the problem of minimizing $f_0:\reals^n\to \Reals$ plus an expectation function defined by $p$ takes the form
\begin{equation}\label{eqn:MinExp}
\nnmin_{x\in\reals^{n}} f_0(x) + \nsum_{i=1}^s p_i f_i(x),
\end{equation}
where $f_i:\reals^n\to \Reals$, $i=1, \dots, s$, represent outcomes occurring with probabilities $p_1, \dots, p_s$. These actual probabilities might be unknown and we examine the effect of replacing them by $p^\nu \in \Delta$, which presumably gets closer to $p$ as $\nu\to \infty$. For example, we may assume initially that $p^\nu_i = 1/s$, because of arguments related to ``maximum entropy,'' but the actual probability vector might have $p_i = 0$ for some $i$ representing outliers or corruption that should not have been in the data set. We may realize this over time and update $p^\nu$ accordingly.

Two-stage stochastic optimization problems take the form \eqref{eqn:MinExp} with $f_0(x)$ being the first-stage cost of decision $x$ and $f_i(x)$ specifying the second-stage cost in scenario $i$, which then occurs with probability $p_i$. The second-stage cost could equal infinity in the absence of complete recourse, which leaves open the possibility that $f_i(x) = \infty$; see for example \cite{chen2022sample} and  \cite[Section 3.J]{primer}. This produces an induced constraint in \eqref{eqn:MinExp} when $p_i >0$. However, if $p_i = 0$ because the underlying data point is an outlier or represents an irrelevant situation, possibly introduced by an adversary, then $p_if_i(x) = 0$ and $x$ might even be a minimizer in \eqref{eqn:MinExp}. An alternative probability $p_i^\nu > 0$, regardless of how close it is to the actual probability $p_i = 0$, would still lead to the conclusion that $x$ is infeasible because $p_i^\nu f_i(x) = \infty$.

Generally, small changes to the probability distribution in \eqref{eqn:MinExp} can shift the solutions significantly when $f_1, \dots, f_s$ are discontinuous or extended real-valued.

\begin{example}{\rm (finite probability space).}\label{eMinExp}
Classification problems with a cross-entropy loss function takes the form \eqref{eqn:MinExp}.  Let $g(x;\xi,\eta)$ be the probability a statistical model (neural net) parameterized by $x\in\reals^n$ assigns to label $\eta\in \{1, 2, \dots, c\}$ given feature vector $\xi\in\reals^d$. For given data $\{(\xi_i, \eta_i )\}_{i=1}^s$, where $\eta_i\in \{1, 2, \dots, c\}$ is the label and $\xi_i\in \reals^d$ is the feature of data point $i$, we obtain a training problem stated as \eqref{eqn:MinExp} with $f_i(x) = -\log g(x;\xi_i,\eta_i)$. In the context of influence functions \cite{koh2017understanding} and outlier analysis, the effect of deviation from the nominal $p_i = 1/s$ becomes important. The case $p_i = 0$ indicates a data point that the training process should ignore and this makes $g(x;\xi_i,\eta_i)=0$ admissible. However, a small change away from $p_i = 0$ would bar any such $x$.

Difficulty also arises when $f_i$ is real-valued but discontinuous such that when \begin{equation}\label{eqn:probexprs}
f_i(x) = H\big( g(\xi_i,x)\big),
\end{equation}
where $H$ is the Heaviside step function, i.e., $H(\gamma) = 1$ if $\gamma >0$ and $H(\gamma) = 0$ if $\gamma \leq 0$, and $g:\reals^m\times\reals^n\to \reals$. This expression occurs in statistical learning with loss functions counting the number of misclassifications or representing the AUC (area under the ROC curve) metric \cite{NortonUryasev.19}. It also arises in engineering design with the objective function then representing a failure probability \cite{royset2006optimal}.
\end{example}
\state Detail. For \eqref{eqn:probexprs}, let $\xi_1 = 0$ and $\xi_2 = 1$, each occurring with probability $p_1 = p_2 = 1/2$, let $g(\xi,x) = \xi + x$, and let $f_0(x) = (1/4)(x-1)^2$ for $x\in [0,1]$ and $f_0(x) = \infty$ otherwise. Then, the objective function in \eqref{eqn:MinExp} has the value $3/4$ at $x = 0$, which is the unique minimizer, and value $(1/4)(x-1)^2 + 1$ for $x\in (0,1]$. Now, suppose that the support of the uncertain parameter is changed from $\{0,1\}$ to $\{1/\nu, 1\}$, while the probabilities remain unchanged. The objective function in \eqref{eqn:MinExp} is now $(1/4)(x-1)^2 + 1$ for $x\in [0,1]$, with a unique minimizer at $x = 1$ regardless of $\nu\in \nats$.\eop

The examples illustrate how a small change to a probability distribution in a problem may cause the values of expectation functions to shift up in a disproportional manner. Since the shift is up, it is futile to attempt to address the situations by passing to a ``robust'' formulation involving worst-case distributions as this would only cause expectation functions to shift even higher. As an alternative, we aim to anticipate the shift from the actual but unknown distribution $P$ to an on-hand distribution $P^\nu$ by considering perturbations of $P^\nu$. While the ideal perturbation that would recover $P$ will remain unknown, a joint optimization over the decision vector $x$ and the perturbation allows us to consider a wider set of solutions and mitigate the effect of the shift from $P$ to $P^\nu$.

\section{Rockafellian Relaxation}\label{sec:rock}

We now abstract the elements of the previous section and formalize the framework. For $\phi:\reals^n\to \Reals$, we consider the {\em actual problem}
\begin{equation}\label{eqn:actualproblem}
\nnmin_{x\in\reals^n} \phi(x).
\end{equation}
However, $\phi$ is unavailable for the decision-maker due to some unsettled or corrupted quantity. To handle this situation we proceed in two steps. First, we design a perturbation scheme for the actual problem and thus place it within a family of problems that are in some sense better positioned to absorb any deviation from the actual problem. Second, we analyze the effect of deviations caused by the unsettled or corrupted quantities. We accomplish the first step by {\em designing} a suitable Rockafellian in anticipation of the deviations of the second step. Generally, a function $f:\reals^m\times\reals^n\to \Reals$ is referred to as a {\em Rockafellian}  for \eqref{eqn:actualproblem} if
\[
f(0, x) = \phi(x)~~~\forall x\in\reals^n.
\]
The first argument of $f$, usually denoted by $u$, specifies a perturbation of the actual problem\footnote{We refer to \cite[Section 5.1]{primer} for a slightly broader definition.}.

We adopt the following notation and terminology: For $\psi:\reals^n\to \Reals$ and $\epsilon \in [0,\infty)$, we define $\inf_x \psi(x) = \inf \{\psi(x)~|~x\in\reals^n\}$, $\dom \psi = \{x\in\reals^n~|~\psi(x)<\infty\}$, and $\epsilon\mbox{-}\hspace{-0.06cm}\nargmin_x \psi(x) = \{\bar x\in \dom \psi~|~\psi(\bar x) \leq \inf_x \psi(x) + \epsilon\}$, with $0\mbox{-}\hspace{-0.06cm}\nargmin_x \psi(x)$ simply being written as $\nargmin_x \psi(x)$. For any set $C$, the {\em indicator function} is given as $\iota_C(x) = 0$ if $x\in C$ and $\iota_C(x) = \infty$ otherwise. We recall that $\psi$ is {\em proper} when $\dom \psi \neq \emptyset$ and $\psi(x)>-\infty$ for all $x\in\reals^n$. It is {\em lower semicontinuous} (lsc) if, for each $x\in\reals^n$,  $\nliminf \psi(x^\nu) \geq \psi(x)$ whenever $x^\nu\to x$. The {\em conjugate} $\psi^*: \reals^n \rightarrow \Reals$ of $\psi$ is given by $\psi^*(y) = \sup_{x} \langle x,y \rangle - \psi(x)$. For $f:\reals^m\times\reals^n\to \Reals$, we let the {\em min-value function} be given by
\[
\cV(u) = \ninf_{x} f(u,x).
\]
If $f$ is a Rockafellian for the actual problem \eqref{eqn:actualproblem}, then $\cV(0)$ specifies its minimum value.

\subsection{Exactness}

While the actual problem \eqref{eqn:actualproblem} can be associated with many Rockafellians, we seek to design one that exhibits an exactness property, which is related to conditions for strong duality \cite[Chapter 11]{VaAn}.

\begin{definition}{\rm (exact Rockafellians).}\label{def:exactness}
A Rockafellian $f:\reals^m\times \reals^n\to \Reals$ for (\ref{eqn:actualproblem}) is {\em exact, supported by} $\bar{y} \in \reals^m$, if
\begin{equation*}
\mathcal{V}(u) \geq \mathcal{V}(0) + \langle \bar y, u\rangle ~~~\forall u\in\reals^m.
\end{equation*}
If the inequality holds strictly for all $u\neq 0$, then $f$ is {\em strictly exact}, supported by $\bar y$.
\end{definition}

A trivial design (but still useful as seen in the next section) is to let $f(u,x) = \infty$ for $u\neq 0$ and $f(0,x) = \phi(x)$, which indeed produces an exact Rockafellian supported by any  $\bar y \in \reals^m$. A characterization of exactness is available via conjugacy. All the proofs in this section appear in the appendix.

\begin{proposition}{\rm (exactness from conjugacy).}\label{lemma:: ExactnessGeom}
Let $f:\reals^m\times \reals^n\to \Reals$ be a Rockafellian for (\ref{eqn:actualproblem}). Then $f$ is exact,  supported by $\bar{y} \in \reals^m$, if and only if $\mathcal{V}(0) = - \mathcal{V}^*(\bar{y})$.
\end{proposition}

Since $-\mathcal{V}^*(\bar{y}) = \inf_{u,x} f(u,x) - \langle \bar y, u\rangle$ and the right-hand side involves joint optimization over $u$ and $x$, one has $\cV(0) \geq -\mathcal{V}^*(\bar{y})$ always and the {\em relaxed problem}
\begin{equation}
\label{eqn:relaxedproblem}
\nnmin_{u\in\reals^m, x \in \reals^n} f(u,x) - \langle \bar y, u\rangle
\end{equation}
is indeed a relaxation of the actual problem \eqref{eqn:actualproblem}. This holds for any $\bar y\in \reals^m$. We refer to the construction of \eqref{eqn:relaxedproblem} as {\em Rockafellian relaxation}. As seen in Proposition \ref{lemma:: ExactnessGeom}, exactness amounts to tightness in this relaxation.

Another characterization of exactness emphasizes relations between minimizers.

\begin{proposition}{\rm (minimizer characterization of exactness).}\label{thm:minChar} Suppose $f:\reals^m\times \reals^n\to \Reals$ is a Rockafellian for the problem \eqref{eqn:actualproblem} of minimizing $\phi:\reals^n\to \Reals$.

If $f$ is exact, supported by $\bar{y}$, then
\begin{equation} \label{eqn:: ExactnessImplication}
x^\star \in \nargmin_x \phi(x) ~~\Longrightarrow ~~ (0, x^\star) \in \nargmin_{u,x} f(u,x) - \langle \bar y, u\rangle.
\end{equation}
When $\nargmin_x \phi(x) \neq \emptyset$, \eqref{eqn:: ExactnessImplication} conversely implies exactness. 

If $f$ is strictly exact, supported by $\bar{y}$, then
\begin{equation} \label{eqn:: StrictExactnessImplication}
(u^\star, x^\star) \in \nargmin_{u,x} f(u,x) - \langle \bar y, u\rangle ~~\Longrightarrow ~~ u^\star = 0, ~x^\star \in \nargmin_x \phi(x).
\end{equation}
When $\nargmin_{u,x} f(u,x) - \langle \bar y, u\rangle \neq \emptyset$ and $\argmin_x f(u,x) \neq \emptyset$ for all $u\in \dom \cV$, then (\ref{eqn:: StrictExactnessImplication}) conversely implies strict exactness.
\end{proposition}

Hence, exactness of a Rockafellian $f$ produces an optimality condition for the actual problem \eqref{eqn:actualproblem}:
\begin{equation}\label{eqn:optcondbasic}
 x^\star \in \nargmin_x \phi(x) ~~\Longrightarrow ~~ \exists u^\star\in\reals^m \mbox{ such that } \,(u^\star, x^\star) \in \nargmin_{u,x} f(u,x) - \langle \bar y, u\rangle.
\end{equation}
The optimality condition helps justify the passing from \eqref{eqn:actualproblem} to \eqref{eqn:relaxedproblem} even without strict exactness. Further discussion of optimality conditions is deferred to Section \ref{sec:optcond}.

There are numerous ways of constructing exact Rockafellians, with some possibilities listed below. Section \ref{sec:appl} gives concrete examples in the context of stochastic optimization under distributional perturbations.

\begin{corollary}{\rm (exactness from convexity).}\label{cor:SuffCondExactConvex}
A Rockafellian $f:\reals^m\times \reals^n\to \Reals$ for \eqref{eqn:actualproblem} is exact, supported by $\bar y\in \reals^m$ when the min-value function $\mathcal{V}$ is convex with $\bar y$ as a subgradient at $0$.

If $\mathcal{V}$ is also strictly convex, then the Rockafellian is strictly exact, supported by $\bar y$.
\end{corollary}

Since an inf-projection of a convex Rockafellian is convex \cite[Proposition 1.21]{primer}, it is immediately clear that the convexity requirement in the corollary holds for such Rockafellians. In this setting, the dual problem of maximizing $\inf_{u,x} f(u,x) - \langle y, u\rangle$ over $y$ furnishes $\bar y$ \cite[Theorem 11.39]{VaAn}.

\begin{proposition}{\rm (exactness from calmness).}\label{prop:SuffCondExactCalm}
Suppose that $f:\reals^m\times \reals^n\to \Reals$ is a Rockafellian for \eqref{eqn:actualproblem}. Let $\|\cdot\|$ be an arbitrary norm on $\reals^m$. For any $\theta\in [0,\infty)$, the function given by $f_\theta(u,x) = f(u,x) + \theta \|u\|$
is also a Rockafellian for \eqref{eqn:actualproblem}. If there is $\bar\theta \in [0,\infty)$ such that
\begin{equation}\label{eqn:calmnesscond}
\cV(u) \geq \cV(0) - \bar\theta \|u\| ~~\forall u\in \reals^m,
\end{equation}
then $f_\theta$ is exact, supported by $\bar y = 0$ as long as $\theta\geq \bar \theta$. The exactness is strict if $\theta>\bar\theta$ and $\cV(0) \in \reals$.
\end{proposition}

The requirement \eqref{eqn:calmnesscond} holds for sufficiently large $\bar \theta$ if the min-value function $\cV$ is calm from below at $0$ \cite[Equation 8(12)]{VaAn} and it is also minorized by some affine function. Calmness in the context of right-hand side perturbations of constraint systems is discussed in \cite{burke1991exact}.

\subsection{Limiting Exactness}

The motivation for designing a Rockafellian $f$ is to mitigate the effect of any unsettled or corrupted quantities in the actual problem \eqref{eqn:actualproblem}. We formalize the setting by considering the approximating functions $f^\nu:\reals^m\times \reals^n\to \Reals$, $\nu \in \nats$. For example, the actual problem and its Rockafellian $f$ might be defined in terms of the ``true'' distribution $P$, while the approximating function $f^\nu$ involves the corrupted distribution $P^\nu$. The approximating functions might also reflect other changes from the Rockafellian, possibly motivated by computational considerations.

These functions lead to the {\em approximating problem}
\begin{equation}\label{eqn:approxproblem}
\nnmin_{u\in\reals^m, x\in\reals^n} f^\nu(u,x) - \langle y^\nu, u\rangle,
\end{equation}
where $y^\nu\in \reals^m$ serves as a substitute for some vector $\bar y$ that in turn supports exactness. In contrast to the actual problem \eqref{eqn:actualproblem} and a relaxed problem \eqref{eqn:relaxedproblem}, which are conceptual in the sense that they involve unsettled quantities, an approximating problem is assumed accessible to the decision-maker and eventually will be solved numerically. Under exactness, there is equivalence between the actual problem and a relaxed problem. We next turn to the  relationship between the chosen relaxed problem \eqref{eqn:relaxedproblem} and its approximations \eqref{eqn:approxproblem} using the notion of epi-convergence\footnote{For background information about epi-convergence; see for example \cite[Chapter 4]{primer}.}.

The functions $\psi^\nu:\reals^n\to \Reals$ are said to {\em epi-converge} to $\psi:\reals^n\to \Reals$ when
\begin{align}
  &\forall x^\nu\to x, ~\nliminf \psi^\nu(x^\nu) \geq \psi(x)\label{eqn:liminfcondition}\\
  &\forall x, ~\exists x^\nu\to x \mbox{ with } \nlimsup \psi^\nu(x^\nu) \leq \psi(x).\label{eqn:limsupcondition}
\end{align}
In the context of Rockafellians, this leads to the following definition.\\

\begin{definition}{\rm (limit-exact Rockafellians).}\label{def:asymRock}
For the problem \eqref{eqn:actualproblem}, the functions $\{f^\nu:\reals^m\times\reals^n\to \Reals, \nu\in\nats\}$ are {\em  limit-exact Rockafellians} if they epi-converge to an exact Rockafellian $f:\reals^m\times\reals^n\to \Reals$ of the problem.

The functions $\{f^\nu, \nu\in\nats\}$ are  {\em strictly limit-exact Rockafellians} if $f$ is strictly exact. The (strict) limiting exactness is said to be supported by $\bar y\in\reals^m$ when $f$ is (strictly) exact, supported by $\bar y$.
\end{definition}

In the presence of limit-exact Rockafellians, we are on solid ground for passing from the actual problem \eqref{eqn:actualproblem} to the approximating problem \eqref{eqn:approxproblem}; the latter provides arbitrarily accurate solutions of the former as seen by the following fact.

\begin{theorem}{\rm (convergence under limiting exactness)}. \label{thm:asymptExact}
For the problem  \eqref{eqn:actualproblem} of minimizing $\phi:\reals^n\to \Reals$, suppose that $\dom \phi \neq \emptyset$ and that $\{f^\nu:\reals^m\times\reals^n\to \Reals, \nu\in\nats\}$ are strictly limit-exact Rockafellians, supported by $\bar y\in \reals^m$. Let $y^\nu\to \bar y$, $\epsilon^\nu\downto 0$, and
\[
(u^\nu,x^\nu) \in \epsilon^\nu\hspace{-0.02cm}\mbox{-}\hspace{-0.06cm}\nargmin_{u,x} \big\{f^\nu(u,x) - \langle y^\nu, u \rangle\big\}.
\]
Then, every cluster point $(\hat u, \hat x)$ of $\{(u^\nu,x^\nu), \nu\in\nats\}$ satisfies $\hat x\in \nargmin_x \phi(x)$.

If strict limiting exactness is replaced by limiting exactness, then $(\hat u, \hat x)$ satisfies the necessary optimality condition \eqref{eqn:optcondbasic} with $f:\reals^m\times\reals^n\to \Reals$ being furnished by the limit of $\{f^\nu, \nu\in\nats\}$.
\end{theorem}

The theorem establishes that the approximating problems \eqref{eqn:approxproblem}, when defined by strictly limit-exact Rockafellians, are in the limit equivalent to the actual problem \eqref{eqn:actualproblem}. They exhibit the desirable continuity property. An approximating problem becomes an accessible and viable alternative to solve in lieu of the unknown actual problem. Still, the approximating problems could be challenging computationally as discussed in Sections 5 and 6.

\section{Applications to Distributional Perturbations}\label{sec:appl}

The wide framework of Section \ref{sec:rock} offers a pathway to analyzing stochastic optimization problems under changes to underlying probability distributions. This section illustrates some possibilities in the context of the ``classical'' stochastic optimization problem \eqref{eqn:MinExp} with perturbations of the probability vector $p$. We supplement the treatment by also considering changes to the support, expectation constraints, and more general distributions.

\subsection{Minimization of Expectations}\label{subsec:MinExp}

In this subsection, we consider the actual problem \eqref{eqn:MinExp} under the assumption that $p\in \Delta$ and each $f_i:\reals^n\to \Reals$ is proper and lsc, $i = 0, 1, \dots, s$. As illustrated in Section \ref{section:: Examples}, a change to $p$ could have outsized effects on solutions. Even if $p^\nu\in \Delta$ is near $p$, we cannot automatically assume that minimizers of the alternative problem
\begin{equation}\label{eqn:naive}
\nnmin_{x\in\reals^n} f_0(x) + \nsum_{i=1}^s p_i^\nu f_i(x)
\end{equation}
would be close to those of \eqref{eqn:MinExp}. It turns out that an approximating problem defined via a Rockafellian is better behaved.

We design the Rockafellian $f:\reals^{s}\times \reals^n\to \Reals$ given by
\begin{equation}\label{eqn:rockExpect1}
f(u,x) = f_0(x) + \nsum_{i=1}^s (p_i + u_i) f_i(x) + \iota_{\{0\}^s}(u).
\end{equation}
For $p^\nu \in \Delta$ and $\theta^\nu \in [0,\infty)$, let
\begin{equation}\label{eqn:rockExpect1nu}
f^\nu(u,x) = f_0(x) + \nsum_{i=1}^s (p_i^\nu + u_i) f_i(x) + \half\theta^\nu\|u\|_2^2 + \iota_{\Delta}(p^\nu + u),
\end{equation}
which, for any $y^\nu\in \reals^s$, defines an approximating problem
\begin{equation}\label{eqn:approxExp}
\nnmin_{u\in\reals^s, x\in \reals^n} f_0(x) + \nsum_{i=1}^s (p_i^\nu + u_i) f_i(x) + \half\theta^\nu\|u\|_2^2 + \iota_{\Delta}(p^\nu + u) - \langle y^\nu, u\rangle.
\end{equation}
In the following, we provide justifications for considering \eqref{eqn:approxExp} as a substitute for the actual problem \eqref{eqn:MinExp}. These approximating problems turn out to have better properties than \eqref{eqn:naive}, which can be sensitive to changes in the probability vector.

\begin{proposition}{\rm (exactness in expectation minimization).}\label{pExactExp}
In the setting of this subsection, with $p\in \Delta$ and proper lsc $f_i:\reals^n\to \Reals$, $i = 0, 1, \dots, s$, the Rockafellian $f$ in \eqref{eqn:rockExpect1} is strictly exact, supported by any $\bar y \in \reals^s$. If $\theta^\nu\to \infty$ and $\theta^\nu \|p^\nu - p\|_2^2 \to 0$ with $p^\nu\in \Delta$, then the functions $\{f^\nu, \nu\in \nats\}$  in \eqref{eqn:rockExpect1nu} are strictly limit-exact Rockafellians with $f$ as their limit.
\end{proposition}
\state Proof. The strict exactness follows immediately from Definition \ref{def:exactness}. The second claim is a special case of Theorem \ref{thm:ExpBoth} below.\eop

The proposition recommends that $\theta^\nu$ should grow slower than $\|p^\nu - p\|_2^2$ vanishes.  The rate by which $p^\nu$ tends to $p$ depends on the setting. In a conceptual study, rate details may not be needed. In algorithmic developments, however, it becomes more critical as the approximating problem \eqref{eqn:approxExp} depends on $\theta^\nu$. One would like to know the magnitude of any corruption or modeling error, but the exact value can, importantly, remain unknown. Approximating probability vectors stemming from empirical distributions based on $\nu$ sample points are supported by the following fact, which stipulates that $\theta^\nu$ should grow slightly slower than $\nu$, but this is just one possibility.

\begin{example}{\rm (convergence of empirical distributions).}\label{prob:Empirical}
Suppose that $\mathbf{p}^\nu$ is the random vector representing the empirical distribution obtained by sampling independently $\nu$ times according to $p \in \Delta$. Then, for any $\epsilon\in (0,1/2)$, one has $\nu^{1/2-\epsilon}\|\mathbf{p}^\nu - p\|_2\to 0$ almost surely.
\end{example}
\state Detail. This fact holds because $p$ can be viewed as the probability mass function of a random variable support on $\{1, 2, \dots, s\}$. Sampling $\nu$ times according to $p$ produces another $s$-dimensional vector specifying the frequency of the different outcomes. We view this vector as a random vector and denote it by ${\bf p}^\nu$. Recall that $\int_1^\infty \exp(-2t^\alpha) dt < \infty$  if $\alpha > 0$ and $\infty$ if $\alpha = 0$. Let $\epsilon \in (0, 1/2)$, $F:\reals\to [0,1]$ be the distribution function corresponding to $p$, and $\mathbf{F}^\nu$ be the (random) distribution function corresponding to $\mathbf{p}^\nu$. By the (univariate) Dvoretzky-Kiefer-Wolfowitz inequality,
\[
\prob\left(\nsup_{t\in\reals} \big|\mathbf{F}^\nu(t) - F(t)\big|  > \frac{\gamma}{\sqrt{\nu}}\right) \leq 2 \exp(-2\gamma^2)\;\;\; \forall \gamma\in\reals.
\]
Setting $\gamma_\nu = \nu^{\epsilon}$, we conclude that $\sum_{\nu=1}^\infty 2 \exp(-2 \gamma_\nu^2) < \infty$ by the integral test for series. The Borel-Cantelli lemma yields then, with probability one, that $\sup_{t\in \reals}|\mathbf{F}^\nu(t) - F(t)| \leq \nu^{-(1/2 - \epsilon)}$ holds for all but finitely many $\nu$. In particular, for any $i$, $|\mathbf{p}_i^\nu - p_i| \leq 2 \nu^{-(1/2 - \epsilon)}$ holds for all but finitely many $\nu$, from which we reach the conclusion.\eop

The minimization over $u$ in the approximating problem \eqref{eqn:approxExp} can largely be carried out ``explicitly,'' which results in a reformulation of the approximating problem with an interpretation as a regularized version of the alternative problem \eqref{eqn:naive}. Specifically, for fixed $x\in\cap_{i=1}^s \dom f_i$, let $F(x) = (f_1(x), \dots, f_s(x))$. Then, the minimization over $u$ in \eqref{eqn:approxExp} amounts to
\begin{align*}
& \inf_{u\in\reals^s} \Big\{ f_0(x) + \nsum_{i=1}^s (p_i^\nu + u_i) f_i(x) + \half\theta^\nu\|u\|_2^2 + \iota_{\Delta}(p^\nu + u) - \langle y^\nu, u\rangle\Big\}\\
& = f_0(x) + \nsum_{i=1}^s p_i^\nu f_i(x) - \sup_{u\in\reals^s} \Big\{ \big\langle y^\nu-F(x), u\big\rangle - \half\theta^\nu\|u\|_2^2 - \iota_{\Delta}(p^\nu + u)\Big\}\\
 &  = f_0(x) + \nsum_{i=1}^s p_i^\nu f_i(x) - r^\nu(x),
\end{align*}
where
\[
r^\nu(x) = \min_{w\in\reals^s} \Big\{ \max_{i=1, \dots, s} w_i - \langle p^\nu, w\rangle + \frac{1}{2\theta^\nu} \big\|y^\nu - F(x) - w\big\|_2^2 \Big\}.
\]
The transition from the sup-expression to $r^\nu(x)$ is achieved by recognizing that the former is the conjugate of the function $u\mapsto \half\theta^\nu\|u\|_2^2 + \iota_{\Delta}(p^\nu + u)$, evaluated at a certain point. This in turn is expressed by the conjugates of $\half\theta^\nu\|\cdot\|_2^2$ and $\iota_{\Delta}$; see \cite[Example 5.29]{primer} as well as Equation 11(3) and Theorem 11.23(a) in \cite{VaAn}.

At least when $f_1, \dots, f_s$ are real-valued, the derivation establishes that the approximating problem \eqref{eqn:approxExp} is equivalently stated as
\[
\nnmin_{x\in\reals^n} f_0(x) + \nsum_{i=1}^s p_i^\nu f_i(x) - r^\nu(x),
\]
where the role of $r^\nu$ as a regularizer is now apparent when comparing to the alternative problem \eqref{eqn:naive}. Since $r^\nu(x) \geq 0$ for all $x\in\reals^n$, the additional term amounts to a ``negative regularization'' that brings down the potentially too high values of the objective function under $p^\nu$.

The value $r^\nu(x)$ equals a Moreau envelope of the real-valued convex function $w\mapsto \max_{i=1, \dots, s} w_i - \langle p^\nu, w\rangle$ evaluated at the point $y^\nu - F(x)$. Thus, by \cite[Theorem 2.26]{VaAn}, $r^\nu$ is continuously differentiable with $\nabla r^\nu(x) = - \nabla F(x)^\top ( y^\nu - F(x) - \hat w(x))/\theta^\nu$, where $\hat w(x)$ is the unique minimizer in the problem defining $r^\nu(x)$, provided that $f_1, \dots, f_s$ are continuously differentiable. Albeit our motivation stems from problems without such differentiability, this derivation helps us understand the characteristics of the regularizer.

Further justification for the passing from the actual problem \eqref{eqn:MinExp} to the approximating problem \eqref{eqn:approxExp} is furnished by the next result on rate of convergence as $p^\nu$ tends to $p$. We adopt the notation $\ball(\bar x,\rho) = \{x\in \reals^n~|~\|x-\bar x\|_2 \leq \rho\}$.  The {\em point-to-set distance} $\dist(\bar x,C) = \inf\{\|x-\bar x\|_2~|~x\in C\}$ if $C\neq\emptyset$ and $\dist(\bar x,C) = \infty$ otherwise.

\begin{theorem}{\rm (rates in expectation minimization).}\label{prop:Rate}
Consider the setting of this subsection with $p,p^\nu\in \Delta$ and proper lsc $f_i:\reals^n\to \Reals$, $i = 0, 1, \dots, s$.  For $\rho\in [0,\infty)$ and $\epsilon\in [0,2\rho]$, suppose that the actual problem \eqref{eqn:MinExp} has a minimizer in $\ball(0,\rho)$ and each one of the approximating problems \eqref{eqn:approxExp} has a minimizer in $\ball(0,\rho)\times\ball(0,\rho)$ with all these problems having minimum values in $[-\rho, \rho-\epsilon]$. If $\{x^\nu, \nu\in\nats\}$ are constructed by
\[
(u^\nu,x^\nu) \in \epsilon\mbox{-}\hspace{-0.06cm}\nargmin_{u,x} f^\nu(u,x) - \langle y^\nu, u\rangle
\]
with $\|u^\nu\|_2 \leq \rho$, $\|x^\nu\|_2 \leq \rho$, and bounded $\{y^\nu,\nu\in\nats\}$, then there are positive constants $\sigma$ and $\tau$ such that for  sufficiently large $\theta^\nu$ and sufficiently small $\|p^\nu-p\|_2$:
\begin{equation}\label{eqn:rateMainRes}
\dist\Big(x^\nu, ~(\epsilon + 2 \eta^\nu)\mbox{-}\hspace{-0.06cm}\nargmin_x f_0(x) + \nsum_{i=1}^s p_i f_i(x)\Big) \leq \eta^\nu,
\end{equation}
where $\eta^\nu = \sigma\|p^\nu-p\|_2 + \max\{\half \theta^\nu\|p^\nu-p\|_2^2,~ \tau/\sqrt{\theta^\nu}\}$. Specifically, one can take $\tau = \beta \sigma$ and
\[
\sigma = \max\Big\{1, ~\nsup_{\nu\in\nats} \|y^\nu\|_2 +  \sqrt{s}\Big(\max\Big\{\kappa, \sqrt{\tfrac{3}{2\alpha}}\beta\Big\} + \kappa\Big)\Big\},
\]
where $\alpha = \min_{i=1, \dots, s} \{p_i~|~p_i > 0\}$, $\beta = \sqrt{2\rho + 2\rho\nsup_{\nu\in\nats}\|y^\nu\|_2 + 4\kappa}$, and $\kappa \in [0,\infty)$ satisfies
\begin{align*}
\ninf\big\{ f_0(x)~\big|~\|x\|_2 \leq \rho\big\} & \geq - \kappa\\
\ninf\big\{ f_i(x)~\big|~\|x\|_2\leq \rho, f_0(x)<\infty, i=1, \dots, s\big\}&\geq -\kappa.
\end{align*}
In terms of these quantities, \eqref{eqn:rateMainRes} holds for $\theta^\nu \geq 9 \beta^2/\alpha^2$ and $\|p^\nu-p\|_\infty \leq \alpha/3$.
\end{theorem}
\state Proof. Suppose that $\nu$ is such that $\theta^\nu \geq 9 \beta^2/\alpha^2$ and $\|p^\nu-p\|_\infty \leq \alpha/3$. Let $\tilde f^\nu(u,x) = f^\nu(u,x) - \langle y^\nu, u\rangle$. The two functions $f$ and $\tilde f^\nu$ are lsc. This allows us to slightly refine \cite[Theorem 6.56]{primer} (details omitted) and confirm that
\begin{equation}\label{eqn:proofRate1}
\dist_*\big((u^\nu,x^\nu), ~\delta\mbox{-}\hspace{-0.06cm}\nargmin_{u,x} f(u,x)\big) \leq \hatsetd_\rho(\epi f, \epi \tilde f^\nu)
\end{equation}
provided that $\delta \geq \epsilon + 2\hatsetd_\rho(\epi f, \epi \tilde f^\nu)$, where $\dist_*(a,A)$ is the point-to-set distance between a point $a \in \reals^s\times\reals^n$ and a set $A\subset\reals^s\times\reals^n$ under the norm $((u,x),(\tilde u, \tilde x)) \mapsto \max\{\|u-\tilde u\|_2, \|x-\tilde x\|_2\}$ and $\hatsetd_\rho(A, B)$ is the truncated Hausdorff distance between sets $A$ and $B$ computed under the norm
\[
\big((u,x,\gamma),(\tilde u, \tilde x,\tilde\gamma) \big) \mapsto\max\big\{\|u-\tilde u\|_2, \|x-\tilde x\|_2, |\gamma - \tilde\gamma|\big\}
\]
on $\reals^s\times\reals^n\times\reals$; see \cite[Section 6.J]{primer} for definitions. Moreover, $\epi f = \{(u,x,\gamma)\in\reals^{s+n+1}~|~f(u,x) \leq \gamma\}$ with a parallel definition for $\epi \tilde f^\nu$.  Since
\[
(\tilde u, \tilde x) \in \delta\mbox{-}\hspace{-0.06cm}\nargmin_{u,x} f(u,x) ~~\Longleftrightarrow~~ \tilde  u = 0, ~\tilde x \in \delta\mbox{-}\hspace{-0.06cm}\nargmin_{x} f_0(x) + \nsum_{i=1}^s p_i f_i(x)
\]
for finite $\delta$, we establish from \eqref{eqn:proofRate1} that
\[
\dist\Big(x^\nu, ~\delta\mbox{-}\hspace{-0.06cm}\nargmin_{x} f_0(x) + \nsum_{i=1}^s p_i f_i(x)\Big) \leq \hatsetd_\rho(\epi f, \epi \tilde f^\nu)
\]
provided that $\infty>\delta \geq \epsilon + 2\hatsetd_\rho(\epi f, \epi \tilde f^\nu)$.

We leverage the Kenmochi condition 6.58 in \cite{primer} to bound $\hatsetd_\rho(\epi f, \epi \tilde f^\nu)$ and proceed in two steps. First, suppose that $(u,x)\in\reals^s\times\reals^n$ satisfies $\|u\|_2 \leq \rho$, $\|x\|_2 \leq \rho$, and $f(u,x) \leq \rho$. Then, $u=0$.
Set
\[
\eta_1^\nu = \max\big\{\half \theta^\nu \|p^\nu-p\|_2^2 + \|y^\nu\|_2 \|p^\nu-p\|_2, ~\|p^\nu-p\|_2\big\}.
\]
It now follows that
\begin{align*}
&\ninf_{\tilde u,\tilde x} \big\{\tilde f^\nu(\tilde u,\tilde x)\,\big|\,\|\tilde u-u\|_2\leq \eta_1^\nu, \|\tilde x-x\|_2 \leq \eta_1^\nu\big\} \leq \tilde f^\nu(p-p^\nu,x)\\
&  ~~~~= f_0(x) + \nsum_{i=1}^s (p_i^\nu + p_i - p_i^\nu) f_i(x) + \half \theta^\nu \|p^\nu-p\|_2^2 - \langle y^\nu, p-p^\nu\rangle\\
& ~~~~\leq f_0(x) + \nsum_{i=1}^s p_i f_i(x)  + \eta_1^\nu = f(u,x) + \eta_1^\nu.
\end{align*}

Second, we derive a similar expression with the roles of $f$ and $\tilde f^\nu$ reversed. Suppose that $(u,x)\in\reals^s\times\reals^n$ satisfies $\|u\|_2 \leq \rho$, $\|x\|_2 \leq \rho$, and $\tilde f^\nu(u,x) \leq \rho$. Then, $x\in \dom f_0$, $p^\nu + u \in \Delta$, and
\[
-2\kappa + \half \theta^\nu\|u\|_2^2 - \rho\|y^\nu\|_2  \leq f_0(x) + \nsum_{i=1}^s (p_i^\nu + u_i) f_i(x) + \half\theta^\nu\|u\|_2^2 + \iota_{\Delta}(p^\nu + u) - \langle y^\nu, u\rangle \leq \rho.
\]
Thus, one obtains $\theta^\nu\|u\|_2^2   \leq 2\rho + 2\rho\|y^\nu\|_2 + 4\kappa$, which means that $\|u\|_2 \leq \beta/\sqrt{\theta^\nu}$.

Next, set $I_0 = \{i~|~p_i =0\}$, $I^+ = \{i~|~p_i > 0\}$, and $I^+_x = \{i \in I^+~|~f_i(x) > 0\}$. Since $|p_i-p_i^\nu| \leq \alpha/3$, $|u_i|\leq \beta/\sqrt{\theta^\nu}$, and $\theta^\nu \geq 9 \beta^2/\alpha^2$, we obtain that $p_i^\nu + u_i = p_i - p_i + p_i^\nu + u_i \geq \alpha - \alpha/3 - \alpha/3 = \alpha/3$ for $i\in I^+$. This implies that for $i^* \in I_x^+$, one has
\begin{align*}
-2\kappa + \tfrac{1}{3}\alpha f_{i^*}(x) - \rho \|y^\nu\|_2
& \leq -2\kappa + \nsum_{i \in I_x^+} (p_i^\nu + u_i) f_i(x) - \rho\|y^\nu\|_2\\
& \leq  f_0(x) + \nsum_{i=1}^s (p_i^\nu + u_i) f_i(x) + \half\theta^\nu\|u\|_2^2 + \iota_{\Delta}(p^\nu + u) - \langle y^\nu, u\rangle \leq \rho.
\end{align*}
Consequently, for all $i\in I^+$, one has $-\kappa \leq f_i(x) \leq \gamma^\nu = 3(\rho + \rho \|y^\nu\|_2  + 2\kappa)/\alpha$. This results in the following bound:
\begin{align*}
\nsum_{i=1}^s p_i f_i(x) - \nsum_{i=1}^s (p_i^\nu + u_i) f_i(x) & =  \nsum_{i \in I^+} (p_i - p_i^\nu - u_i) f_i(x) - \nsum_{i \in I_0} (p_i^\nu + u_i) f_i(x)\\
& \leq \max\{\kappa,\gamma^\nu\}\nsum_{i \in I^+} |p_i - p_i^\nu - u_i| + \kappa\nsum_{i \in I_0} (p_i^\nu + u_i) \\
& \leq \eta_2^\nu = \sqrt{s} \big(\max\{\kappa,\gamma^\nu\} + \kappa\big)\big(\|p^\nu-p\|_2 + \beta/\sqrt{\theta^\nu} \big).
\end{align*}
We leverage this bound to compute
\begin{align*}
&\ninf_{\tilde u, \tilde x} \big\{f(\tilde u, \tilde x)~\big|~\|\tilde u - u\|_2\leq \beta/\sqrt{\theta^\nu}, \|\tilde x - x\|_2 \leq \beta/\sqrt{\theta^\nu} \big\} \leq f(0,x) = f_0(x) + \nsum_{i=1}^s p_i f_i(x)\\
 & = \tilde f^\nu(u,x) + \nsum_{i=1}^s p_i f_i(x) - \nsum_{i=1}^s (p_i^\nu + u_i) f_i(x) - \half \theta^\nu \|u\|_2^2 + \langle y^\nu, u\rangle\\
& \leq \tilde f^\nu(u,x) + \eta_2^\nu + \beta\|y^\nu\|_2 /\sqrt{\theta^\nu}.
\end{align*}
The Kenmochi condition 6.58 in \cite{primer} then produces the bound
\[
\hatsetd_\rho(\epi f, \epi \tilde f^\nu) \leq \max\big\{\eta_1^\nu, ~\beta/\sqrt{\theta^\nu}, ~\eta_2^\nu + \beta\|y^\nu\|_2 /\sqrt{\theta^\nu}\big\}.
\]
The conclusion follows after slightly relaxing the upper bound to facilitate a simple expression and noting that $\gamma^\nu \leq \beta\sqrt{3/(2\alpha)}$.\eop

The trade-offs in the rate expression in the theorem could be addressed, for example, by setting
$\theta^\nu = \|p^\nu-p\|_2^{-4/3}$. Then, $1/\sqrt{\theta^\nu}$ and $\theta^\nu \|p^\nu-p\|_2^2$ decay at the same rate as $\|p^\nu-p\|_2^{2/3}$. We conclude that the error, in the sense of the theorem, by passing from the actual problem \eqref{eqn:MinExp} to the approximating problem \eqref{eqn:approxExp} vanishes at the rate $\|p^\nu-p\|_2^{2/3}$. This holds under mild assumptions; $f_0, f_1, \dots, f_s$ could be nonconvex, nonsmooth, discontinuous, and extended real-valued. In comparison, the alternative problem \eqref{eqn:naive} may exhibit large errors in the absence of additional assumptions.

\subsection{Minimization of Expectations: Refinements}\label{subsec:refine}

We now discuss three refinements of the previous subsection: (A) approximating functions $f^\nu$ based on $\Phi$-divergence instead of $\|\cdot\|_2^2$, (B) perturbation of probability vector {\em and} corresponding support, and (C) strictly limit-exact Rockafellians without relying on the rate of decay of $\|p^\nu - p\|_2$.\\

\noindent {\bf Refinement A.} While retaining the Rockafellian $f$ from Subsection \ref{subsec:MinExp}, we consider new approximating functions $f^\nu$. We recall that the $\Phi$-divergence of $q\in \Delta$ from $\bar q\in \Delta$ is given by
\[
d_\Phi(q|\bar q) = \nsum_{i=1}^s \bar q_i \Phi(q_i/\bar q_i),
\]
where $\Phi:\reals\to \Reals$ is a convex function, real-valued on $[0,\infty)$, with a unique minimizer at $1$ and a minimum value of $0$. (While $\Phi$-divergence can be defined for slightly broader classes of convex functions $\Phi$, we limit the scope to such functions here.) In the expression for $\Phi$-divergence, we adopt the usual conventions: $0 \cdot \Phi(0/0) = 0$, and, for $\beta \in (0,\infty)$, $0 \cdot \Phi(\beta/0) = \beta \nlim_{\alpha\to \infty} \Phi(\alpha)/\alpha$. The requirements on $\Phi$ are satisfied by the commonly used Kullback-Leibler divergence, Burg-entropy, $J$-divergence, $\chi^2$-distance, modified $\chi^2$-distance, variational distance, and Hellinger distance.

As an alternative to the path in Subsection \ref{subsec:MinExp}, we now define
\begin{equation}\label{eqn:phifnu}
f^\nu(u, x) = f_0(x) + \nsum_{i=1}^s (p^\nu_i + u_i) f_i(x) + \theta^\nu d_\Phi(p^\nu + u|p^\nu) + \iota_{\Delta}(p^\nu + u).
\end{equation}
The resulting approximating problem is also considered in \cite{gotoh2021data} for a study of out-of-sample performance, but there under the assumptions of strict convexity and twice smoothness of each $f_i$. We require only proper lsc functions to confirm that the new approximating functions are strictly limit-exact Rockafellians.

\begin{theorem}{\rm (exactness in expectation minimization; $\Phi$-divergence).}
In the setting of Subsection \ref{subsec:MinExp}   with $p,p^\nu\in \Delta$ and proper lsc $f_i:\reals^n\to \Reals$, $i = 0, 1, \dots, s$, consider $f$ from \eqref{eqn:rockExpect1} and $f^\nu$ from \eqref{eqn:phifnu}. The functions $\{f^\nu, \nu\in \nats\}$ are strictly limit-exact Rockafellians with $f$ as their limit provided that $\theta^\nu\to \infty$ and $\theta^\nu \|p^\nu - p\|_2 \to 0$.
\end{theorem}
\state Proof. We show that $f^\nu$ epi-converges to $f$ via the conditions \eqref{eqn:liminfcondition} and  \eqref{eqn:limsupcondition}.

Let $(u^\nu, x^\nu) \to (u,x)$. First, suppose that $u\neq 0$. Then, $f(u,x) = \infty$. We consider two cases. (i) Suppose that $p + u \not\in \Delta$. Then, $p^\nu + u^\nu \not\in \Delta$ for sufficiently large $\nu$ because $p^\nu + u^\nu \to p + u$. Thus, $\nliminf f^\nu(u^\nu,x^\nu) = \infty$. (ii) Suppose that $p + u \in \Delta$. We assume without loss of generality that $p^\nu + u^\nu \in \Delta$ for all $\nu$ because $f^\nu(u^\nu, x^\nu)=\infty$ when $p^\nu + u^\nu \not\in \Delta$. Next, there is $i^*\in \{1, \dots, s\}$ such that $p_{i^*}>0$ and $u_{i^*} \neq 0$. This claim holds by the following argument:

For the sake of contradiction, suppose that $u_i = 0$ for all $i$ with $p_i>0$. Since $p+u\in \Delta$, one has
\[
\nsum_{i=1}^s p_i + u_i = \nsum_{i|p_i>0} p_i  + \nsum_{i|p_i=0} u_i = 1.
\]
This implies that $\nsum_{i|p_i=0} u_i =0$. All summands in this sum cannot be zero because $u\neq 0$. In turn, this means that there is a negative summand, say $j$, which makes $p_j + u_j = 0 + u_j < 0$. However, this contradicts the fact that $p+u \in \Delta$.

Returning to the main argument, the claim confirms that $(p_{i^*}^\nu + u_{i^*}^\nu)/p_{i^*}^\nu \to 1 + u_{i^*}/p_{i^*} \neq 1$. The continuity of $\Phi$ on $(0,\infty)$, the convexity of $\Phi$, and the unique minimizer of $\Phi$ at 1 with minimum value of 0 imply that
\[
\Phi\Big(\frac{p_{i^*}^\nu + u_{i^*}^\nu}{p_{i^*}^\nu}\Big) \to \beta > 0,
\]
where $\beta = \Phi(1 + u_{i^*}/p_{i^*})$ if $1 + u_{i^*}/p_{i^*}>0$ and $\beta = \lim_{\alpha\downto 0} \Phi(\alpha)$ otherwise. Thus,
\[
\theta^\nu p_{i^*}^\nu \Phi\Big(\frac{p_{i^*}^\nu + u_{i^*}^\nu}{p_{i^*}^\nu}\Big) \to \infty.
\]
Since $\Phi(\alpha) \geq 0$ for all $\alpha\in\reals$, $\theta^\nu d_\Phi(p^\nu + u^\nu|p^\nu)\to \infty$. This fact together with the assumption about proper lsc functions $f_0, f_1, \dots, f_s$ imply that
\begin{align*}
&\nliminf f^\nu(u^\nu, x^\nu)\\
 & \geq \nliminf f_0(x^\nu) + \nsum_{i=1}^s \nliminf \big((p^\nu_i + u_i^\nu) f_i(x^\nu)\big) + \nliminf \big( \theta^\nu d_\Phi(p^\nu + u^\nu|p^\nu)\big) + \nliminf\iota_{\Delta}(p^\nu + u^\nu)\\
& \geq f_0(x) + \nsum_{i=1}^s (p_i + u_i) f_i(x) + \nliminf \big( \theta^\nu d_\Phi(p^\nu + u^\nu|p^\nu)\big) = \infty.
\end{align*}

Second, suppose that $u = 0$. Similar to above, one obtains
\begin{align*}
\nliminf f^\nu(u^\nu, x^\nu) & \geq f_0(x) + \nsum_{i=1}^s p_i f_i(x) + \nliminf \big( \theta^\nu d_\Phi(p^\nu + u^\nu|p^\nu)\big)\\
& \geq f_0(x) + \nsum_{i=1}^s p_i f_i(x) = f(0,x).
\end{align*}
We have shown that the liminf-condition \eqref{eqn:liminfcondition} holds for $f^\nu$ and $f$.

Next, we turn to the limsup-condition \eqref{eqn:limsupcondition}. Let $(u,x) \in \reals^s\times \reals^n$.  If $u\neq 0$, then $f(u,x) = \infty$ and the limsup-condition holds trivially. Thus, we concentrate on the case $u = 0$. Set $x^\nu = x$ and $u^\nu = p - p^\nu$. Then,
\[
f^\nu(u^\nu, x^\nu) = f^\nu(p - p^\nu, x) = f(0,x) + \theta^\nu d_\Phi(p|p^\nu).
\]
Let $i\in \{1, \dots, s\}$ be arbitrary. We consider two cases: (i) Suppose that $p_i > 0$. Then, for sufficiently large $\nu$, $p_i^\nu > 0$. Since $\Phi$ is convex and real-valued in a neighborhood of 1 and $\Phi(1) = 0$, there exist $\epsilon, \delta\in (0,\infty)$ such that $\Phi(\alpha) \leq \delta|\alpha-1|$ for all $\alpha \in [1-\epsilon , 1+\epsilon]$. This implies that, for sufficiently large $\nu$, one has
\[
p_i^\nu \Phi\Big(\frac{p_i}{p_i^\nu}\Big) \leq p_i^\nu \delta \Big| \frac{p_i}{p_i^\nu} - 1 \Big| = \delta | p_i -  p_i^\nu |.
\]
(ii) Suppose that $p_i = 0$. If $p_i^\nu > 0$, then $p_i^\nu \Phi(p_i/p_i^\nu) = p_i^\nu \Phi(0) =  \Phi(0) |p_i - p_i^\nu|$. If $p_i^\nu = 0$, then, by convention, $p_i^\nu \Phi(p_i/p_i^\nu) = 0 \cdot \Phi(0/0) = 0$. Since $i$ is arbitrary and $\theta^\nu \|p^\nu - p\|_2\to 0$, we conclude that $\theta^\nu d_\Phi(p^\nu|p) \to 0$. This means that $\nlimsup f^\nu(u^\nu,x^\nu) \leq f(u,x)$. We have established that $f^\nu$ epi-converges to $f$.\eop\\

\noindent {\bf Refinement B.} Subsection \ref{subsec:MinExp} deals with changes to the probability vector $p$ in the actual problem \eqref{eqn:MinExp}. Now, we consider the possibility that the support may also be unsettled. Changes to a support underpin adversarial training of neural nets \cite{Madry.18}, but with a focus on conservativeness. We consider relaxations.

For proper lsc functions $f_0:\reals^n\to \Reals$ and $g:\reals^m \times \reals^n\to \Reals$, support $\{\xi_1, \dots, \xi_s\} \subset \reals^m$, and corresponding probability vector $p = (p_1, \dots, p_s)\in \Delta$, we consider the actual problem
\begin{equation}\label{eqn:actualExpSup}
\nnmin_{x\in \reals^n} f_0(x) + \nsum_{i=1}^s p_i g(\xi_i,x).
\end{equation}
The problem furnishes additional structure to \eqref{eqn:MinExp}, where one has $f_i(x) = g(\xi_i,x)$. Mimicking the development in Subsection \ref{subsec:MinExp}, we design the Rockafellian $f:\reals^{s + sm}\times \reals^n\to \Reals$ given by
\begin{equation}\label{eqn:fsupport}
f\big((u,v),x\big) = f_0(x) + \nsum_{i=1}^s (p_i + u_i) g(  \xi_i + v_i,x) + \iota_{\{0\}^s}(u) + \iota_{\{0\}^{sm}}(v).
\end{equation}
For $p^\nu = (p_1^\nu, \dots, p_s^\nu)\in \Delta$, $\xi^\nu = (\xi_1^\nu, \dots, \xi_s^\nu)\in \reals^{sm}$, and $\theta^\nu,\lambda^\nu \in [0,\infty)$, we define the approximating functions by
\begin{equation}\label{eqn:fnusupport}
f^\nu\big((u,v),x\big) = f_0(x) + \nsum_{i=1}^s (p_i^\nu + u_i) g(\xi_i^\nu  + v_i,x) + \half\theta^\nu\|u\|_2^2 + \half\lambda^\nu\|v\|_2^2 + \iota_{\Delta}(p^\nu + u).
\end{equation}
These approximating functions turn out to be  strictly limit-exact Rockafellians.

\begin{theorem}{\rm (exactness in expectation minimization; support uncertainty).}\label{thm:ExpBoth}
Suppose that $f_0:\reals^n\to \Reals$ and $g:\reals^m \times \reals^n\to \Reals$ are proper lsc functions and $p,p^\nu\in \Delta$. For the actual problem \eqref{eqn:actualExpSup}, the Rockafellian $f$ in \eqref{eqn:fsupport} is strictly exact, supported by any $\bar y \in \reals^{s+sm}$.

If $\theta^\nu,\lambda^\nu\to \infty$, $\theta^\nu \|p^\nu - p\|_2^2 \to 0$, and $\lambda^\nu \|\xi^\nu -   \xi\|_2^2 \to 0$, where $\xi = (\xi_1, \dots, \xi_s)\in \reals^{sm}$, then the functions $\{f^\nu, \nu\in \nats\}$ in \eqref{eqn:fnusupport} are strictly limit-exact Rockafellians with $f$ as their limit.
\end{theorem}
\state Proof. The claim about strict exactness holds trivially. We establish that $f^\nu$ epi-converges to $f$ using \eqref{eqn:liminfcondition} and \eqref{eqn:limsupcondition}. Let $(u^\nu,v^\nu,x^\nu)\to (u,v,x)$. Then, $\nliminf f_0(x^\nu)\geq f(x)>-\infty$, $\nliminf \half\theta^\nu\|u^\nu\|_2^2 \geq \iota_{\{0\}^s}(u)$, $\nliminf \half\lambda^\nu\|v^\nu\|_2^2 \geq \iota_{\{0\}^{sm}}(v)$, and  $\nliminf \iota_{\Delta}(p^\nu + u^\nu ) \geq \iota_{\Delta}(p + u )$ because $\Delta$ is closed. Since $f^\nu((u^\nu,v^\nu),x^\nu) = \infty$ if $p^\nu + u^\nu \not\in \Delta$, we assume without loss of generality that $p^\nu + u^\nu \in \Delta$ for all $\nu\in\nats$. This in turn implies that
\[
\nliminf \nsum_{i=1}^s (p_i^\nu + u_i^\nu) g(\xi_i^\nu  + v_i^\nu,x^\nu) \geq \nsum_{i=1}^s \nliminf  (p_i^\nu + u_i^\nu) g(\xi_i^\nu  + v_i^\nu,x^\nu) \geq \nsum_{i=1}^s (p_i + u_i) g(  \xi_i  + v_i,x).
\]
We combine these facts to conclude that the liminf-condition \eqref{eqn:liminfcondition}  holds. For the limsup-condition \eqref{eqn:limsupcondition}, fix $(u,v,x)$ and construct $x^\nu = x$, $u^\nu = p - p^\nu$, and $v^\nu =   \xi - \xi^\nu$. Without loss of generality, we assume that $u = 0$ and $v = 0$ because $f((u,v),x)= \infty$ otherwise. Since $p\in \Delta$, $p^\nu + u^\nu \in \Delta$. Then,
\begin{align*}
&\nlimsup f^\nu\big((u^\nu,v^\nu),x^\nu\big)\\
& \leq f_0(x) + \nlimsup \nsum_{i=1}^s (p_i^\nu + u_i^\nu) g(\xi_i^\nu  + v_i^\nu,x) + \nlimsup \half\theta^\nu\|u^\nu\|_2^2 + \nlimsup \half\lambda^\nu\|v^\nu\|_2^2\\
&= f_0(x) + \nsum_{i=1}^s p_i g(  \xi_i,x) + \nlimsup \half\theta^\nu\|p- p^\nu\|_2^2 + \nlimsup \half\lambda^\nu\|  \xi - \xi^\nu\|_2^2\\
& = f_0(x) + \nsum_{i=1}^s p_i g(  \xi_i,x) = f\big((u,v),x\big).
\end{align*}
We conclude that $f^\nu$ epi-converges to $f$ and the assertion follows by Definition \ref{def:asymRock}.\eop\\

\noindent {\bf Refinement C.} We return to the actual problem \eqref{eqn:MinExp}, but now design an alternative Rockafellian with the advantage that in implementation of the resulting approximating problem there is no need to know the rate by which $p^\nu$ tends to $p$. Specifically, for $\theta \in [0,\infty)$, we define a Rockafellian given by
\begin{equation}\label{eqn:rockRefiC}
f(u,x) = f_0(x) + \nsum_{i=1}^s (p_i + u_i) f_i(x) + \theta \|u\|_1  + \iota_\Delta(p+u).
\end{equation}
In contrast to earlier Rockafellians defined in terms of indicator functions, it is not immediately clear whether $f$ is exact. We also define the approximating functions by setting
\begin{equation}\label{eqn:arockRefiC}
f^\nu(u,x) = f_0(x) + \nsum_{i=1}^s (p_i^\nu + u_i) f_i(x) + \theta \|u\|_1  + \iota_\Delta(p^\nu+u).
\end{equation}
We see next that the desired exactness emerges when $\theta$ is sufficiently large.

\begin{theorem}{\rm (exactness under alternative Rockafellian).}\label{tAltRock}
In the setting of Subsection \ref{subsec:MinExp} with $p,p^\nu\in \Delta$ and proper lsc $f_i:\reals^n\to \Reals$, $i = 0, 1, \dots, s$, consider $f$ from \eqref{eqn:rockRefiC} and $f^\nu$ from \eqref{eqn:arockRefiC}. Suppose that there are $\eta\in [0,\infty)$ and $\bar x\in\reals^n$ such that $f_i(\bar x) <\infty$  and $\ninf_x f_i(x) \geq -\eta$, $i = 0, 1, \dots, s$. 

If $\theta$ is sufficiently large, then the Rockafellian $f$ is strictly exact, supported by $0\in \reals^s$, and, under the additional assumption that $ \|p^\nu - p\|_2 \to 0$, the approximating functions $\{f^\nu, \nu\in \nats\}$ are strictly limit-exact Rockafellians with $f$ as their limit.
\end{theorem}
\state Proof. First, we consider the strict exactness of $f$. For $u\in \reals^s$, we adopt the notation
\[
\delta(u) = \ninf_x f_0(x) + \nsum_{i=1}^s (p_i + u_i) f_i(x) + \iota_\Delta(p+u).
\]
If for some $\kappa \in [0,\infty)$ one has $\delta(u) \geq \delta(0) - \kappa \|u\|_1$ for all $u\in\reals^s$ in a neighborhood of $0$ and $\delta(0)\in\reals$, then one can invoke Proposition \ref{prop:SuffCondExactCalm} to conclude that $f$ is strictly exact, supported by $0 \in \reals^s$, provided that $\theta$ is sufficiently large. Certainly,
\[
-2\eta \leq \delta(0) \leq f_0(\bar x) + \nsum_{i=1}^s p_i f_i(\bar x)<\infty.
\]
For the purpose of constructing $\kappa$, we adopt the notation
$I_0 = \{i~|~p_i =0\}$, $I^+ = \{i~|~p_i > 0\}$, and $\alpha = \min\{p_i~|~i\in I^+\}$.

Suppose that $\|u\|_1 \leq \alpha/2$. If $\delta(u) \geq \delta(0)$, then $\kappa = 0$ applies so we proceed under the assumption that $\delta(u) < \delta(0)$. Since $\delta(0)$ is finite, this implies that $p+u\in \Delta$. Let $\epsilon \in (0,\delta(0) - \delta(u)]$. Since $-2\eta \leq \delta(u) < \delta(0)<\infty$, there exists $\hat x\in\reals^n$ such that
\[
f_0(\hat x) + \nsum_{i=1}^s (p_i + u_i) f_i(\hat x) \leq \delta(u) + \epsilon \leq \delta(0).
\]
Let $i\in I^+$. Then, one has
\[
  (p_i + u_i) f_i(\hat x) \leq \delta(0) - f_0(\hat x) - \nsum_{j\neq i} (p_j + u_j) f_j(\hat x) \leq \delta(0) + 2\eta.
\]
The probability $p_i + u_i\geq \alpha/2$ because $p_i\geq \alpha$ and $|u_i| \leq \alpha/2$. Consequently,
\[
f_i(\hat x) \leq \frac{\delta(0) + 2\eta}{p_i + u_i} \leq \frac{2\delta(0) + 4\eta}{\alpha}<\infty.
\]
This bound together with $f_i(\hat x)\geq -\eta$ facilitate the following development:
\begin{align*}
  \delta(u) &\geq f_0(\hat x) + \nsum_{i=1}^s (p_i + u_i) f_i(\hat x) - \epsilon\\
   & = f_0(\hat x) + \nsum_{i \in I^+} (p_i + u_i) f_i(\hat x) + \nsum_{i \in I_0} (p_i + u_i) f_i(\hat x) - \epsilon\\
   & \geq f_0(\hat x) + \nsum_{i \in I^+} p_i f_i(\hat x) + \nsum_{i \in I^+}  u_i f_i(\hat x) - \eta \nsum_{i \in I_0} (0 + u_i) - \epsilon\\
&\geq f_0(\hat x) + \nsum_{i \in I^+} p_i f_i(\hat x) - \max\Big\{\eta, \frac{2\delta(0) + 4\eta}{\alpha}\Big\} \|u\|_1 - \eta \|u\|_1 - \epsilon\\
&\geq \ninf_x f_0(x) + \nsum_{i=1}^s p_i f_i(x) - \bigg(\max\Big\{\eta, \frac{2\delta(0) + 4\eta}{\alpha}\Big\} + \eta\bigg) \|u\|_1 - \epsilon.
\end{align*}
Since $\epsilon$ is arbitrary, we have shown that $\delta(u) \geq \delta(0) - \kappa \|u\|_1$ when $\|u\|_1 \leq \alpha/2$ for $\kappa$ given inside the parentheses on the previous line.

Second, we consider the approximating functions and establish epi-convergence using the liminf-condition \eqref{eqn:liminfcondition} and the limsup-condition \eqref{eqn:limsupcondition}. Let $\theta\in [0,\infty)$. For the liminf-condition, let $(u^\nu, x^\nu)\to (u,x)$. If $p+u \not\in\Delta$, then $p^\nu + u^\nu \not\in\Delta$ for sufficiently large $\nu$ because $\Delta$ is closed. This implies that
\begin{equation}\label{eqn:liminfProofAlt}
\nliminf f^\nu(u^\nu,x^\nu)  \geq  f(u,x)
\end{equation}
holds because each side equals infinity.  If $p+ u\in \Delta$, then, without loss of generality, we assume that $p^\nu + u^\nu\in \Delta$ for all $\nu$. We obtain that
\begin{align*}
\nliminf f^\nu(u^\nu, x^\nu) & \geq  \nliminf f_0(x^\nu) + \nsum_{i=1}^s \nliminf (p_i^\nu + u_i^\nu) f_i(x^\nu) + \nliminf \theta \|u^\nu\|_1 + \nliminf \iota_\Delta(p^\nu + u^\nu)\\
& \geq  f_0(x) + \nsum_{i=1}^s \nliminf (p_i^\nu + u_i^\nu) f_i(x^\nu) + \theta\|u\|_1.
\end{align*}
Since $p_i^\nu + u_i^\nu\geq 0$, we find that $\nliminf (p_i^\nu + u_i^\nu) f_i(x^\nu) \geq (p_i + u_i) f_i(x)$. Thus, \eqref{eqn:liminfProofAlt} holds again.

For the limsup-condition, let $(u,x)\in \reals^s\times \reals^n$ and construct $x^\nu = x$ and $u^\nu = p - p^\nu + u$. We assume without loss of generality that $p + u\in \Delta$. Then, $p^\nu + u^\nu = p + u \in \Delta$. Moreover,
\begin{align*}
\nlimsup f^\nu(u^\nu,x^\nu) & \leq  f_0(x) + \nsum_{i=1}^s (p_i + u_i) f_i(x) + \nlimsup \theta \|p-p^\nu + u\|_1\\
& \leq  f_0(x) + \nsum_{i=1}^s (p_i + u_i) f_i(x) + \nlimsup \theta \|p^\nu-p\|_1  + \theta \|u\|_1\\
& =  f_0(x) + \nsum_{i=1}^s (p_i + u_i) f_i(x) + \theta \|u\|_1 = f(u,x).
\end{align*}
We conclude that $f^\nu$ epi-converges to $f$ regardless of $\theta$.\eop

\begin{example}{\rm (bankrupt-prone newsvendor).} In response to a random demand, a newsvendor orders $x$ newspaper each morning at a cost of $\gamma>0$ dollars per paper. The newsvendor sells a newspaper for $\delta>\gamma$ dollars. Unsold papers cannot be returned and are worthless at the end of the day. The goal of the newsvendor is to choose an order quantity that minimizes the expected loss (expense minus income). Classically, the  newsvendor can handle a loss of any size (see, e.g., \cite[Section 1.C]{primer}), but we adjust the problem by requiring the loss to be almost surely nonpositive and show that Rockafellian relaxation leads to a decision that remains hidden for an analyst following a more naive approach.
\end{example}
\state Detail. When ordering $x$ newspapers and $\xi$ is the demand, the loss turns out to be
\[
\begin{cases}
   \gamma x -\delta x    &  \mbox{ if } x \leq \xi\\
   \gamma x - \delta \xi &  \mbox{ if } \xi < x \leq \delta \xi/\gamma\\
   \infty                &  \mbox{ otherwise.}
\end{cases}
\]
The threshold $\delta \xi/\gamma$ is the largest order quantity that avoids the expense $\gamma x$ to exceed the income $\delta \xi$.

We assume that the demand has a finite distribution with support $\{\xi_i \in \nats, i = 1, \dots, s\}$ and associated probabilities $\{p_i, i = 1, \dots, s\}$. The problem of minimizing the expected loss can then be written as the following two-stage stochastic optimization problem (cf. \cite[Section 1.C]{primer})
\begin{equation*}\label{newsmodel}
\nnmin_{x\in\nats} \; (\gamma-\delta) x + \nsum_{i=1}^s p_i f_i(x),
\end{equation*}
where the recourse function is $f_i(x) = \inf\{\delta y^\lplus \; | \; y^\lplus -y^\lminus = x - \xi_i, \; (\delta -\gamma)x/\delta \geq y^\lplus\geq 0,\; y^\lminus\geq 0\}$. For the specific instance with $\gamma = 1$, $\delta = 2$, $\{\xi_i = i, i=1, 2, 3\}$, and $p = (0, 1/4, 3/4)$, we find that the unique minimizer is $x = 3$ with minimum value $-5/2$. However, if we use the slightly different probabilities $p^\nu = (1/\nu, 1/4 - 1/\nu, 3/4)$ for $\nu \geq 9$ in place of $p$, then the unique
minimizer becomes $x = 2$ with minimum value $-2 + 2/\nu$. This holds even as $\nu\to \infty$ and there is a disconcerting sensitivity to arbitrarily small changes to the probability vector. The minimizer $x = 3$ of the actual problem remains hidden for an analyst that  proceeds naively with the seemingly accurate probability vector $p^\nu$.

The Rockafellian given by $f(u,x) = \iota_\nats(x) + (\gamma-\delta) x + \nsum_{i=1}^s (p_i + u_i) f_i(x) + \theta \|u\|_1 + \iota_\Delta(p+u)$ and its approximation $f^\nu$ obtained by replacing $p$ by $p^\nu$ are special cases of \eqref{eqn:rockRefiC} and \eqref{eqn:arockRefiC}. The resulting approximating problem
has minimizer $(u^\nu,x^\nu) = (-1/\nu, 0, 1/\nu, 3)$ with minimum value $-5/2 + 6/\nu$ for $\nu\geq 9$ and $\theta \geq 2$. Thus, the approximating problem allows us to recover the minimizer of the actual problem despite having incorrect probabilities.\eop

\subsection{Expectations in Compositions}\label{subsec:constraints}

Expectation functions also arise in constraints due to fairness considers in statistical learning (see Example \ref{eFairSVM}), in the construction of Neyman-Pearson classifiers \cite{rigollet2011neyman}, and in reliability constraints involving buffered failure probabilities \cite{RockafellarRoyset.10}. Again, small changes to the underlying probabilities may cause large errors in solutions of the corresponding optimization problems. In this subsection, we consider the actual problem
\begin{equation}\label{eqn:svm}
\nnmin_{x\in\reals^n} f_0(x) + h\Big(\nsum_{i=1}^s p_i G_{i}(x)\Big),
\end{equation}
where $p\in \Delta$, $f_0:\reals^n\to \Reals$ and $h:\reals^m\to \Reals$ are proper lsc functions, and $G_i:\reals^n\to \reals^m$ is continuous for each $i$, with $G_i(x) = (g_{i1}(x), \dots, g_{im}(x))$ and $g_{ik}:\reals^n\to \reals$. We are especially interested in $h = \iota_{(-\infty,0]^m}$, which
produces the expectation constraints
\[
\nsum_{i=1}^s p_i g_{ik}(x) \leq 0, ~~k=1, \dots, m,
\]
but the results below hold for any proper lsc $h$. In contrast to earlier sections, we limit the treatment to real-valued expectation functions to focus on the restrictions imposed via $h$.

Since perturbation of a constraint can cause large changes in minimizers of the corresponding problems, it is clear that naively replacing $p$ by $p^\nu\in \Delta$ in \eqref{eqn:svm} may result in large errors. We again turn to a Rockafellian and design
\begin{equation}\label{eqn:Rockcomp1}
f(u,x) = f_0(x) + h\Big(u + \nsum_{i=1}^s p_i G_{i}(x) \Big) + \iota_{\{0\}^m}(u).
\end{equation}
For $p^\nu\in \Delta$ and $\theta^\nu \in [0,\infty)$, we select the approximating functions given by
\begin{equation}\label{eqn:Rockcomp1nu}
f^\nu(u,x) =   f_0(x) + h\Big(u + \nsum_{i=1}^s p_i^\nu G_{i}(x)\Big) + \half\theta^\nu \|u\|_2^2.
\end{equation}
Again strict limiting exactness holds under mild assumptions.

\begin{proposition}{\rm (exactness in composite optimization).}\label{prop:exactComposite}
In the setting of this subsection with $p,p^\nu\in \Delta$, proper lsc functions $f_0:\reals^n\to \Reals$ and $h:\reals^m\to \Reals$, and continuous mappings $G_i:\reals^n\to \reals^m$, $i=1, \dots, s$, consider $f$ from \eqref{eqn:Rockcomp1} and $f^\nu$ from \eqref{eqn:Rockcomp1nu}. Then, the Rockafellian $f$ is strictly exact, supported by any $\bar y \in \reals^m$. If $\theta^\nu\to \infty$ and $\theta^\nu \|p^\nu - p\|_2^2\to 0$, then the functions $\{f^\nu, \nu\in\nats\}$ are strictly limit-exact Rockafellians with $f$ as their limit.
\end{proposition}
\state Proof. The strict exactness follows directly from Definition \ref{def:exactness}. For the second assertion, we leverage the liminf-condition \eqref{eqn:liminfcondition} and the limsup-condition \eqref{eqn:limsupcondition} to establish that $f^\nu$ epi-converges to $f$. Let $(u^\nu,x^\nu) \to (u,x)\in\reals^m\times \reals^n$. Then, $\nliminf f_0(x^\nu) \geq f_0(x)$, $\nliminf \half \theta^\nu \|u^\nu\|_2^2 \geq \iota_{\{0\}^m}(u)$, and
\[
u^\nu + \nsum_{i=1}^s p_i^\nu G_{i}(x^\nu) \to u+  \nsum_{i=1}^s p_i G_{i}(x).
\]
Since $f_0$ and $h$ are proper, this implies that
\[
\nliminf \bigg( f_0(x^\nu) + h\Big(u^\nu + \nsum_{i=1}^s p_i^\nu G_{i}(x^\nu)\Big) + \half\theta^\nu \|u^\nu\|_2^2\bigg) \geq f_0(x) + h\Big(u + \nsum_{i=1}^s p_i G_{i}(x) \Big) + \iota_{\{0\}^m}(u).
\]
For $(u,x)\in\reals^m\times \reals^n$, we construct $x^\nu = x$ and
\[
u^\nu =  u + \nsum_{i=1}^s p_i G_{i}(x) -\nsum_{i=1}^s p_i^\nu G_{i}(x).
\]
If $u \neq 0$, then $\nlimsup f^\nu(u^\nu,x^\nu) \leq f(u,x)$. If $u = 0$, then
\[
\nlimsup f^\nu(u^\nu,x^\nu) \leq f_0(x) + h\Big(\nsum_{i=1}^s p_i G_{i}(x) \Big) + \nlimsup \half\theta^\nu \Big\|\nsum_{i=1}^s p_i G_{i}(x) -\nsum_{i=1}^s p_i^\nu G_{i}(x)\Big\|_2^2.
\]
Since $\theta^\nu \|p^\nu - p\|_2^2 \to 0$, the last term on the right-hand side vanishes and we have established that $\nlimsup f^\nu(u^\nu,x^\nu) \leq f(u,x)$. Thus, $f^\nu$ epi-converges to $f$.\eop

We observe that the result holds without any constraint qualification. Thus, it addresses difficult cases with $h = \iota_{(-\infty,0]^m}$ and $\nsum_{i=1}^s p_i G_{i}(x) \geq 0$ for all $x\in\reals^n$.

If $h$ is not only proper and lsc but also convex, then the minimization over $u$ in the approximating problem can be achieved ``explicitly.'' (The case $h = \iota_{(-\infty, 0]^m}$ is a prime example.) Specifically, for any $x\in\reals^n$ and $v\in\reals^m$, one has
\begin{align*}
\ninf_{u} \big\{ f_0(x) + h(u + v) + \half\theta^\nu \|u\|_2^2 - \langle y^\nu, u\rangle\big\} &= f_0(x) - \nsup_{u} \big\{ \langle y^\nu, u\rangle - h(u + v) - \half\theta^\nu \|u\|_2^2\big\}\\
& = f_0(x) + h^\nu(v),
\end{align*}
where
\[
h^\nu(v) = - \nmin_{w\in\reals^m} \Big\{ h^*(w) - \langle v, w\rangle + \frac{1}{2\theta^\nu} \|y^\nu - w\|_2^2 \Big\}.
\]
The transition from the sup-expression to $h^\nu(x)$ is achieved by recognizing that the former is the conjugate of the function $u\mapsto  h(u +v ) + \half\theta^\nu\|u\|_2^2$ evaluated at $y^\nu$. This in turn is expressed by the conjugate of $h$ and the conjugate of $\half\theta^\nu\|\cdot\|_2^2$; see \cite[Example 5.29]{primer} as well as 11(3) and Theorem 11.23(a) in \cite{VaAn}.
Thus, for any $y^\nu\in \reals^m$, the approximating problem
\[
\nnmin_{u\in\reals^m, x\in\reals^n} f_0(x) + h\Big(u + \nsum_{i=1}^s p_i^\nu G_{i}(x)\Big) + \half\theta^\nu \|u\|_2^2 - \langle y^\nu, u\rangle
\]
is equivalently expressed as
\[
\nnmin_{x\in\reals^n} f_0(x) + h^\nu \Big(\nsum_{i=1}^s p_i^\nu G_{i}(x)\Big).
\]
In view of Proposition \ref{prop:exactComposite}, this problem is better behaved than the naive alternative obtained by simply replacing $p$ by $p^\nu$ in the actual problem \eqref{eqn:svm}. Thus, a change of $p$ in \eqref{eqn:svm} should also be accompanied by a change in $h$, to the better behaved $h^\nu$. The latter is in fact continuously differentiable with
gradient
\[
\nabla h^\nu(v) = \hat w, ~\mbox{ where } \hat w \in \nargmin_{w\in\reals^m} \Big\{ h^*(w) - \langle v, w\rangle + \frac{1}{2\theta^\nu} \|y^\nu - w\|_2^2 \Big\}.
\]
This holds by the inversion rule for subgradients \cite[Proposition 5.37]{primer}.

While evaluation of $h^\nu$ requires the solution of a convex problem, it reduces to quadratic programming when $h = \iota_{(-\infty,0]^m}$. Then, $h^\nu(v) = \max_{w \geq 0} \langle v, w\rangle - \frac{1}{2\theta^\nu} \|y^\nu - w\|_2^2$.

The approximating function $f^\nu$ in this subsection is closely related to augmented Lagrangians within the general duality framework of \cite[Chapter 11]{VaAn}. This observation also points to the possibility of developing exact Rockafellians by replacing $\iota_{\{0\}^m}(u)$ with $\theta\|u\|_1$ in the definition of $f$. Via Proposition \ref{prop:SuffCondExactCalm}, this leads to exactness results parallel to that in Theorem \ref{tAltRock}, at least when each $G_i$ as well as $f_0$ are continuously differentiable and the ``standard'' qualification \cite[Equation (4.16)]{primer} holds.

\subsection{Countable Supports}

As an illustration of possibilities beyond finite probability spaces, we consider countable supports and thus address Example \ref{eConvDistr}. For proper lsc functions $f_i: \reals^n\to \Reals$, $i= 0, 1, 2, \dots$, consider the actual problem
\begin{equation}\label{eqn:actualCount}
\nnmin_{x\in\reals^n} f_0(x) + \nsum_{i=1}^\infty p_i f_i(x).
\end{equation}
Let $\ell^1 = \{(u_1, u_2, \dots)~|~\sum_{i=1}^\infty |u_i|<\infty\}$ and adopt the norm given by $\|u\|_{\ell^1} = \sum_{i=1}^\infty |u_i|$. With slight abuse of notation, we let $\Delta = \{q \in \ell^1~|~\|q\|_{\ell^1} = 1, q_i \geq 0, i\in \nats\}$. Suppose that $p \in \Delta$.

While our original definition of Rockafellians only allows for finite-dimensional perturbations, it trivially extends to perturbations defined on $\ell^1$. In particular, the function $f:\ell^1\times \reals^n\to \Reals$ given by
\[
f(u,x) = f_0(x) + \nsum_{i=1}^\infty (p_i+u_i) f_i(x) + \iota_{\{0\}}(u)
\]
can be viewed as a Rockafellian for the actual problem \eqref{eqn:actualCount}. Although we do not formally define the framework, it is clear that minimizing the actual problem is equivalent to minimizing this Rockafellian. Thus, we have, in an extended sense, strict exactness.

We know from Example \ref{eConvDistr} that naively replacing $p$ by some approximating probabilities $p^\nu \in \Delta$ in the actual problem \eqref{eqn:actualCount} may cause significant errors even if $\|p^\nu - p\|_{\ell^1}\to 0$. It turns out that the Rockafellian supports the development of better behaved approximating problems.

The probability vector $p^\nu$ defines the approximating function given by
\[
f^\nu(u,x) = f_0(x) +  \nsum_{i=1}^\infty (p_i^\nu+u_i) f_i(x) + \theta^\nu \|u\|_{\ell^1} + \iota_\Delta(p^\nu + u).
\]
In view of the strict exactness (in an extended sense) of $f$, a justification for passing from the actual problem \eqref{eqn:actualCount} to minimizing $f^\nu$ hinges on epi-convergence of $f^\nu$ to $f$. The convergence in the liminf- and limsup-conditions \eqref{eqn:liminfcondition} and \eqref{eqn:limsupcondition} are now understood in the sense of the norm $(u,x)\mapsto \max\{\|u\|_{\ell^1}, \|x\|_2\}$. Under such epi-convergence, cluster points of minimizers of $f^\nu$ would specify a minimizer of the actual problem; see \cite{Royset.18} for details about epi-convergence.

\begin{proposition}{\rm (epi-convergence under countable support).}\label{pCountSupport}
In the setting of this subsection with $p,p^\nu\in \Delta$ and proper lsc functions $f_i: \reals^n\to \Reals$, $i= 0, 1, 2, \dots$, suppose that $\theta^\nu \to \infty$, $\theta^\nu\|p^\nu - p\|_{\ell^1}\to 0$, and, for each $\bar x\in \reals^n$, there exist $\epsilon>0$ and $\kappa \in \mathbb{R}$ such that $f_i(x) \geq \kappa$ for all $i \in \nats$ and $x$ with $\|x-\bar x\|_2 \leq \epsilon$.
Then, $f^\nu$ epi-converges to $f$.
\end{proposition}
\state Proof. Towards establishing the liminf-condition \eqref{eqn:liminfcondition}, suppose that $\|u^\nu - u\|_{\ell^1}\to 0$ and $\|x^\nu - x\|_2 \to 0$.
Since $f^\nu(u^\nu,x^\nu) = \infty$ when $p^\nu + u^\nu \not\in\Delta$, we assume without loss of generality that $p^\nu + u^\nu \in \Delta$. We consider two cases. First, if $u \neq 0$, then $\|u^\nu\|_{\ell^1} \geq \half\|u\|_{\ell^1}>0$ for sufficiently large $\nu$. Thus, $\theta^\nu \|u^\nu\|_{\ell^1} \to \infty$. Moreover, $\nsum_{i=1}^\infty (p_i^\nu + u_i^\nu) f_i(x^\nu) \geq \kappa$ for sufficiently large $\nu$. These facts lead to
\[
\nliminf f^\nu(u^\nu, x^\nu) \geq \nliminf f_0(x^\nu) + \nliminf \Bigg( \nsum_{i=1}^\infty (p_i^\nu + u_i^\nu) f_i(x^\nu) \Bigg) + \nliminf \theta^\nu \|u^\nu\|_{\ell^1} = \infty = f(u,x).
\]
Second, suppose that $u = 0$. We note that $\{p^\nu + u^\nu, \nu\in\nats\}$ can be viewed as a sequence of finite measures on the measurable space\footnote{We denote the collection of all subsets of $\nats$ by $\mathcal{P}(\nats)$.} $(\nats, \mathcal{P}(\nats))$ that converge weakly to $p$. Let $g^\nu,g:\nats \to \Reals$ be defined by $g(i) = f_i(x)$  and $g^\nu(i) = f_i(x^\nu)$ for all $i\in \nats$. We equip $\nats$ with the usual distance $(i,j)\mapsto |i - j|$, which makes it a metric space. Since $f_i$ is lsc,
\[
\liminf_{\substack{\nu\to \infty \\ |j-i|\to 0}}  g^\nu(j) \geq g(i).
\]
Thus, by the monotonicity of expectations (see for example \cite[Proposition 8.53]{primer}), one has
\begin{equation}\label{eqn:countliminf}
\nsum_{i=1}^\infty p_i  \liminf _{\substack{\nu\to \infty \\ |j-i|\to 0}}  g^\nu(j) \geq \nsum_{i=1}^\infty p_i g(i).
\end{equation}
We obtain from \cite[Theorem 2.4]{Feinberg_Fatou} that
\begin{equation}\label{eqn:countliminf2}
\nliminf \nsum_{i=1}^\infty (p_i^\nu + u_i^\nu) g^\nu (i) \geq \nsum_{i=1}^\infty p_i  \liminf _{\substack{\nu\to \infty \\ |j-i|\to 0}}  g^\nu(j)
\end{equation}
under the condition:
\[
\lim_{\beta\to \infty} \limsup_{\nu\to \infty} \nsum_{i \in \nats | g_-^\nu(i) \geq \beta} (p_i^\nu + u_i^\nu)g_-^\nu(i) = 0,
\]
where $g_-^\nu(i) = - \min\{0, g^\nu(i)\}$. To see that this condition holds, observe that $\sup_{i \in \nats} g_-^\nu(i) \leq \max\{0, -\kappa\}$ for sufficiently large $\nu$; hence, we find that $\{i\in \nats~|~ g_-^\nu(i) \geq \beta\} = \emptyset$ for sufficiently large $\nu$ when $\beta > \max\{0, -\kappa\}$.
We now combine \eqref{eqn:countliminf} and \eqref{eqn:countliminf2} to reach
\[
\nliminf \nsum_{i=1}^\infty (p_i^\nu + u_i^\nu)f_i(x^\nu) \geq \nsum_{i=1}^\infty p_i f_i(x).
\]
This in turn allows us to conclude that $\nliminf f^\nu(u^\nu,x^\nu) \geq f(u,x)$.

For the limsup-condition \eqref{eqn:limsupcondition}, it suffices to consider $u = 0$ and $x\in\reals^n$. Construct $x^\nu = x$ and $u^\nu = p - p^\nu$. Then, $\|u^\nu\|_{\ell^1} \to 0$ and $p^\nu + u^\nu \in \Delta$. We now have
\begin{align*}
\nlimsup f^\nu(u^\nu,x^\nu) & \leq f_0(x) +  \nlimsup \Big(\nsum_{i=1}^\infty (p_i^\nu+u_i^\nu) f_i(x)\Big) + \nlimsup \big(\theta^\nu \|p^\nu-p\|_{\ell^1}\big)\\
& = f_0(x) +  \nsum_{i=1}^\infty p_i f_i(x) =  f(0,x)
\end{align*}
because $\theta^\nu \|p^\nu-p\|_{\ell^1} \to 0$.\eop

We note that the lower-boundedness by $\kappa$ in the proposition is closely related to equi-lsc  of the class of functions $\{f_i, i\in \nats\}$ \cite[Section 7.C]{VaAn}, which can be relaxed further to asymptotically equi-lsc.

\section{First-Order Optimality Conditions}\label{sec:optcond}

A Rockafellian associated with a problem defines a necessary first-order optimality conditions for the problem under a ``standard'' qualification; see \cite[Theorem 5.10]{primer}. Exactness serves as an alternative qualification.

We adopt the following terminology. For a function $\psi:\reals^n\to \Reals$ and a point $x$ at which it is finite, the set of {\em subgradients} of $\psi$ at $x$ is denoted by $\partial \psi(x)$.  We set $\partial \psi(x) = \emptyset$ for other $x$. For a set $C\subset \reals^n$ and one of its points $x$, we denote by $N_C(x)$ the {\em normal cone} to $C$ at $x$. We set $N_C(x) = \emptyset$ for other $x$. These concepts are understood in the general sense of \cite{VaAn}. The {\em graph} of $\partial \psi$ is $\gph \partial \psi = \{(x,v) \in\reals^n\times \reals^n ~|~v\in \partial \psi(x)\}$.
The {\em outer limit} of a sequence of sets $\{S^\nu\subset \reals^n, \nu\in\nats\}$ is defined\footnote{The outer limit of $S^\nu$ is denoted by $\nlimsup S^\nu$ in \cite{VaAn}.} as
\[
\nOutLim S^\nu = \{x\in \reals^n~|~\exists \mbox{ subsequence } N\subset \nats \mbox{ and $x^\nu \in S^\nu$, $\nu \in N$, such that } x^\nu \Nto x\}.
\]

A first-order optimality condition for the actual problem now follows straightforwardly.

\begin{theorem}{\rm (first-order optimality condition).}\label{thm:firstOptimal} For the problem of minimizing a proper function $\phi:\reals^n\to \Reals$, suppose that $f:\reals^m\times \reals^n\to \Reals$ is an exact Rockafellian supported by $\bar y\in \reals^m$. Then, the following necessary {\em first-order optimality condition} holds:
\[
x^\star \in \nargmin_x \phi(x) ~~\Longrightarrow ~~ \exists u^\star \in \reals^m \mbox{ such that }\,(\bar y, 0) \in \partial f(u^\star,x^\star).
\]
\end{theorem}
\state Proof. Since $x^\star \in \nargmin_x \phi(x)$ implies that
$(0, x^\star) \in \nargmin_{u,x} f(u,x) - \langle \bar y, u\rangle$
by Proposition \ref{thm:minChar}, the Fermat rule (see for example \cite[Theorem 4.73]{primer}) implies, in turn, that $(\bar y, 0) \in \partial f(0,x^\star )$ provided that $f(0, x^\star)\in \reals$. The latter holds because $\phi$ is proper. Thus, $0$ furnishes $u^\star$ in the assertion.\eop

In the earlier sections, we passed from a Rockafellian to its approximations $f^\nu:\reals^m\times\reals^n\to \Reals$, which presumably are minimized in approximating problems. Now, suppose that the approximating problem \eqref{eqn:approxproblem} is only solved in the sense of the optimality condition $(y^\nu, 0) \in \partial f^\nu(u,x)$. Would resulting solutions tend to solutions satisfying the optimality condition $(\bar y, 0) \in \partial f(u,x)$ in Theorem \ref{thm:firstOptimal}?
This might be answered in the affirmative as seen next.

\begin{theorem}{\rm (convergence to first-order optimality condition).}\label{thm:solutionFirstOptimal} For the problem of minimizing a proper function $\phi:\reals^n\to \Reals$, suppose that $y^\nu \to \bar y$, $\epsilon^\nu\downto 0$, and $f:\reals^m\times \reals^n\to \Reals$ is an exact Rockafellian supported by $\bar y\in \reals^m$. If the functions $\{f^\nu:\reals^m\times \reals^n\to \Reals, \nu\in\nats\}$ satisfy the property
\begin{equation}\label{eqn:outerlimitOpt}
\nOutLim (\gph \partial f^\nu) \subset \gph \partial f
\end{equation}
and $(u^\nu,x^\nu)$ solves $(y^\nu, 0) \in \partial f^\nu(u,x)$ with tolerance $\epsilon^\nu$, i.e., $\dist((y^\nu,0), \partial f^\nu(u^\nu,x^\nu)) \leq \epsilon^\nu$, then every cluster point of $\{(u^\nu,x^\nu), \nu\in\nats\}$ satisfies the necessary first-order optimality condition in Theorem \ref{thm:firstOptimal}.

Under the additional assumptions that $\{f^\nu, \nu\in\nats\}$ are strictly limit-exact Rockafellians with $f$ being their limit and each $f^\nu$ being proper, lsc, convex, then \eqref{eqn:outerlimitOpt} holds automatically and every cluster point $(\hat u, \hat x)$ of $\{(u^\nu,x^\nu), \nu\in\nats\}$ satisfies $\hat u = 0$ and $\hat x \in \nargmin_x \phi(x)$.
\end{theorem}
\state Proof. Under \eqref{eqn:outerlimitOpt}, the first conclusion follows immediately. The second conclusion holds by invoking Attouch's theorem; see \cite[Theorem 7.41]{primer}.
It applies because $f^\nu$ being proper, lsc, and convex implies that $f$ is lsc and convex by virtue of being the limit of $f^\nu$. Moreover, $\epi f \neq \emptyset$ because $\phi$ is proper. The possibility $f(u,x) = -\infty$ is ruled about because it would have implied that $f$ is nowhere finite. Thus, $f$ is also proper.
\eop

We end by examining the requirement \eqref{eqn:outerlimitOpt} in the setting of Subsection \ref{subsec:MinExp}. Suppose that $f_1, \dots, f_s$ are continuously differentiable and $\partial f_0$ is outer semicontinuous, which will be the case if $f_0$ is convex or if $f_0$ is an indicator function or in many other case. Then, \eqref{eqn:outerlimitOpt}  holds and Theorem \ref{thm:solutionFirstOptimal} applies. This justifies a computational approach involving the solution of the approximating problem (\ref{eqn:approxproblem}) in the sense of the optimality condition $(y^\nu, 0) \in \partial f^\nu(u,x)$, which might be much more viable than (global) minimization. For large $\nu$, the solution would be close to satisfying the optimality condition in Theorem \ref{thm:firstOptimal} for the actual problem.

The claim about \eqref{eqn:outerlimitOpt} holds by the following arguments. Let $F(x) = (f_1(x), \dots, f_s(x))$. Under a sum rule (see, for example, \cite[Proposition 4.67]{primer}), for $x\in \dom f_0$,
\[
\partial f(0,x) =  \reals^s \times  \Big( \nsum_{i=1}^s p_i \nabla f_i(x) + \partial f_0(x) \Big)
\]
and $\partial f(u,x) = \emptyset$ for $u\neq 0$. Similarly, for $x\in \dom f_0$ and $u \in -p^\nu + \Delta$, one obtains via \cite[Theorem 4.64]{primer} that
\[
\partial f^\nu(u,x) =  \big( F(x) + \theta^\nu u + N_{\Delta}(p^\nu+u) \big) \times  \Big( \nsum_{i=1}^s (p_i^\nu + u_i) \nabla f_i(x) + \partial f_0(x) \Big).
\]
For other $(u,x)$, $\partial f^\nu(u,x) = \emptyset$. Suppose that $(u,x,v,w) \in \nOutLim (\gph \partial f^\nu)$. Then, there exist a subsequence $N\subset \nats$ and points $(u^\nu,x^\nu,v^\nu,w^\nu) \in \gph \partial f^\nu \Nto (u,x,v,w)$. Thus, $v^\nu - F(x^\nu) - \theta^\nu u^\nu \in N_{\Delta}(p^\nu+u^\nu)$  and $w^\nu - \nsum_{i=1}^s (p_i^\nu + u_i^\nu) \nabla f_i(x^\nu) \in \partial f_0(x^\nu)$. Since $\partial f_0$ is outer semicontinuous, one has
\[
w - \nsum_{i=1}^s (p_i + u_i) \nabla f_i(x) \in \partial f_0(x).
\]
Since $N_{\Delta}(p^\nu+u^\nu)$ is a cone, we also have $v^\nu/\theta^\nu - F(x^\nu)/\theta^\nu - u^\nu \in N_{\Delta}(p^\nu+u^\nu)$. The normal cone mapping $N_{\Delta}$ is outer semicontinuous, which implies that $- u \in N_{\Delta}(p+u)$ because $\theta^\nu\to \infty$. Thus, $p+u \in \Delta$ because otherwise $N_\Delta(p+u)$ would have been empty. The set $\Delta$ is nonempty and convex. Thus, $- u \in N_{\Delta}(p+u)$ if and only if $\langle -u, q - (p+u)\rangle \leq 0$ for all $q\in \Delta$. In particular, $\langle -u, p - (p+u)\rangle \leq 0$, which implies that $\|u\|_2 \leq 0$. We have shown that $u = 0$. Thus, $(u,x,v,w) \in \gph \partial f$, and (\ref{eqn:outerlimitOpt}) holds in this case.

\section{Numerical Results}\label{sec:numerical}

The computational challenges of solving an approximating problem \eqref{eqn:approxproblem} depend on the nature of the chosen Rockafellian. In the setting of Subsection \ref{subsec:constraints}, the resulting approximating problem connects with augmented Lagrangians and associated computational approaches become available; see, e.g., \cite[Subsection 6.B]{primer}. For the expectation minimization problem \eqref{eqn:MinExp}, Rockafellians lead to \eqref{eqn:approxExp} and \eqref{eqn:arockRefiC} with the potentially challenging term $(p_i^\nu + u_i) f_i(x)$, which effectively amounts to a bilinear term because it can be replaced by $(p_i^\nu + u_i) z_i$ and the additional constraint $f_i(x)\leq z_i$. This special structure can be exploited in various ways including by means of McCormick relaxations; see, e.g., \cite{ScottStuberBarton.11}. In the context of large-scale problem instances such as from statistical learning, block-coordinate descent algorithms and stochastic proximal-gradient methods \cite{BolteSabachTeboulle.14,AravkinDavis.19,NutiniLaradjiSchmidt.22} are especially promising for solving \eqref{eqn:approxExp} and \eqref{eqn:arockRefiC}; one would cycle between optimizing $u$ and $x$. Optimizing $u$ amounts to linear or quadratic optimization and optimizing $x$ resembles the actual problem for which one can assume there is a suitable algorithm. We utilize this approach in three numerical examples of statistical learning with corruption.

\subsection{Statistical Learning with Corruption}

Empirical risk minimization for statistical learning leads to problems of the form \eqref{eqn:MinExp}, with $s$ being the number of training data points, $f_i(x)$ is the loss for data point $i$ under the statistical model (neural net) represented by $x$, and $f_0(x)$ specifies a regularization term; see Example \ref{eMinExp}. Typically one assigns each data point the probability $1/s$, but this would be ``incorrect'' if some of the data were corrupted. With $C\subset \{1, \dots, s\}$ being the set of corrupted data, the ``correct'' probabilities to assign the data points would be
\[
p = (p_1, \dots, p_s) ~~\mbox{ with } ~~ p_i = \begin{cases}
  0 & \mbox{ if } i \in C\\
  1/(s - |C|) & \mbox{ if } i \in \{1, \dots, s\} \setminus C,
\end{cases}
\]
where $|C|$ is the cardinality of $C$. These probabilities are not available because one cannot easily identify which data points are corrupted. We only have access to the ``incorrect'' probability vector $(1/s, \dots, 1/s)$, which thus plays the role of $p^\nu$ in the notation of Subsection \ref{subsec:MinExp}. The earlier discussion identifies the flaw of simply replacing $p$ by $p^\nu$ in \eqref{eqn:MinExp}, which we below refers to as {\em empirical risk minimization}.
As an alternative, we adopt the Rockafellian of Refinement C in Subsection \ref{subsec:refine} and minimize $f^\nu$ from \eqref{eqn:arockRefiC} with $p^\nu = (1/s, \dots, 1/s)$ using a block-coordinate descent heuristic.\\

\noindent {\bf Heuristic Algorithm for Rockafellian Relaxation.} Minimize $f^\nu$ in \eqref{eqn:arockRefiC} approximately by setting $u^0 = 0$, selecting an initial $x^0$, and assigning $k = 0$. Then, in iteration $k$:
\begin{enumerate}
    \item Carry out $\sigma$ epochs of SGD\footnote{SGD stands for ``stochastic gradient descent'' but this is doubly misleading as it is neither a descent method nor involves gradients only; $f_i$ is nonsmooth and thus we need to consider subgradients \cite[Section 3.G]{primer}.} as applied to the function $x\mapsto \sum_{i=1}^s (p_i^\nu + u_i^k) f_i(x)$ with $x^k$ as initial point. Let $x^{k+1}$ be the point upon termination.

    \item  Compute $u^\star \in \nargmin_{u} \sum_{i=1}^s (p_i^\nu + u_i) f_i(x^{k+1}) + \theta \|u\|_1 + \iota_\Delta(p^\nu + u)$ using the simplex method. Set $u^{k+1} = \mu u^\star + (1-\mu)u^k$, replace $k$ by $k+1$, and go to Step 1.
\end{enumerate}
Below, we use $\sigma = 10$.  The default stepsize parameter $\mu = 0.5$ and the penalty parameter $\theta = 0.4$. (These values are chosen after pilot runs for 12 different pairs of values of $\mu$ and $\theta$.)

\subsection{Experiments}

We experiment on two data sets from computer vision (MNIST \cite{lecun2010mnist} and CIFAR \cite{Krizhevsky09learningmultiple}) and one data set from text analytics (IMDB \cite{maas-EtAl:2011:ACL-HLT2011}).\\

\noindent {\bf MNIST Experiments.} For the MNIST data set, we consider the digits 0, 1, 2 that produce 18623 training data points of which a fraction is assigned corrupted (i.e., incorrect) labels randomly. We adopt a fully connected neural net with three hidden layers consisting of 320, 320, and 200 units, respectively. The hidden layers has the ReLU activation function and the 3-unit output layer utilizes the softmax function. The number of trainable weights in the neural net is 417880, which then is the dimension of $x$. The loss function is cross-entropy as in Example \ref{eMinExp}. The image input consists of 28-by-28 pixels with the pixel values being standardized prior to model input.

Using Tensorflow 2.10.0, the runtime to execute 50 iterations of the heuristic (for a total of 500 epochs) is 800 seconds of which 300 seconds represent the $u$-optimization. (The $u$-portion can be reduced significantly by switching to a subgradient method, but this appears unnecessary when the training data set is of moderate size as here.)

\begin{figure}
\centering
\includegraphics[width=0.99\textwidth]{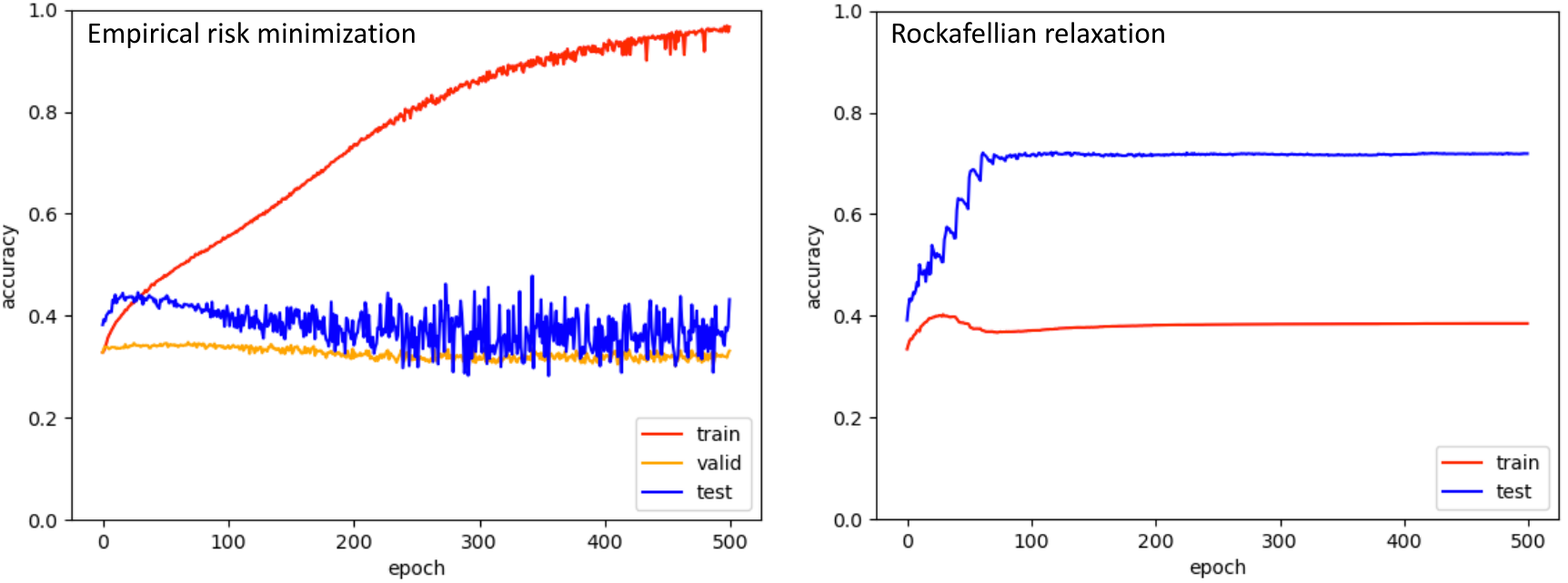}
\caption{Training, validation, and test accuracy for empirical risk minimization (left) and for Rockafellian relaxation (right) on MNIST with 65\% corrupted labels.}\label{mnist}
\end{figure}

Figure \ref{mnist}(right) shows training accuracy (red curve) during the 500 epochs of the heuristic when applied to a training data set with 65\% corruption, i.e., 12104 and 6519 of the training data points are corrupted and clean, respectively. We note that training accuracy is the fraction of the full training data set for which the neural net makes a prediction that matches that of the training data label. (Of course, 65\% of those labels are corrupted.) The blue curve in the figure gives the test accuracy during the 500 epochs using another 3147 clean data points. It plateaus around 0.722. The initial sawtooth shape of the blue curve stems from the $u$-optimization every 10th epoch.

For comparison, Figure \ref{mnist}(left) shows parallel results obtained using empirical risk minimization\footnote{There are numerous approaches to label noise in machine learning such as noisy channels, data pruning, and regularization; see \cite{NortonRoyset.22,chen2022labelpaper,Narasimhanetal.23} and references therein. We omit a comprehensive review and comparison.}, which assigns equal probability ($1/s$) for all data points. The training accuracy (red curve) keeps improving but the test accuracy (blue curve) fails to improve consistently after a few epochs because the neural net is being fitted to the corrupted data points. As is common practice, we allocate 20\% of the training data for validation and this can guide stopping. The orange line in Figure \ref{mnist}(left) reports the validation accuracy. Its peak corresponds to a neural net with test accuracy of 0.421. (The test accuracy peaks at 0.504, but this represents an overly optimistic value that is hard to achieve in practice.) The column labeled ``65\%'' in Table \ref{tab:MNISTaccuracy} reports these values. We conclude that Rockafellian relaxation obtains significantly better test accuracy than empirical risk minimization and largely avoids the need for holding out a portion of the training data set for validation. Figure \ref{mnist}(right) shows that the heuristic used for solving the Rockafellian relaxation problem can terminate whenever the training accuracy stabilizes.

\begin{table}[ht]
\centering
\begin{tabular}{l||r|r|r|r|r|r}
    &  \multicolumn{5}{c}{Percentage corrupted training data}\\
\hline
Method                               & 65\%          & 60\%        & 50\%        & 25\%          & 10\%\\
\hline\hline
Empirical risk minimization           &0.421-0.504	&0.738-0.809  &0.949-0.953	&0.988-0.989	&0.992-0.995\\
Rockafellian ($\mu = 0.5, \theta = 0.4$)  &0.722  	    &0.973	      &0.990        &0.995	        &0.993
\end{tabular}
\caption{Test accuracy for empirical risk minimization and Rockafellian relaxation on MNIST under different levels of corruption.}
\label{tab:MNISTaccuracy}
\end{table}

Table \ref{tab:MNISTaccuracy} summarizes corresponding results under lower levels of corruption. Again, the range of values in the empirical-risk-minimization row specifies test accuracy when the validation accuracy peaks (lower value) and the practically unattainable maximum test accuracy (upper value). Rockafellian relaxation retains a significant advantage for $60\%$ and 50\% corruption, but it vanishes for small levels of corruption as could be expected. In practice, one would not know the level of corruption and it thus seems prudent to adopt Rockafellian relaxation as a means to safeguard against any potential level of corruption.

Returning to the case with 65\% corruption, Table \ref{tab:u} reports the evolution of the $u$-vector across its 49 updates. The columns labeled ``1. iteration'' show the distribution of $u_i$-values after the first $u$-optimization across the 12104 corrupted data points and 6519 clean data points. A portion of the corrupted data points are assigned $u_i = -2.7\cdot 10^{-5}$, while a smaller number of the clean data points are assigned the same value. The columns labeled ``2. iteration'' show similar numbers after the second update of the $u$-vector. The columns labeled ``49. iteration'' display the final values of $u_i$. A large portion ($10197/12104$) of the corropted data points has $u_i = -5.4\cdot 10^{-5}$, which means that they have been completely removed from consideration because this value cancels the nominal probability $p_i^\nu = 5.4\cdot 10^{-5}$. A small fraction of the clean data points receive the same treatment, but the vast majority ($5134/6519$) remain with their nominal probability, or nearly so. This explains the performance of Rockafellian relaxation in Table \ref{tab:MNISTaccuracy}: the approach automatically identifies and ``removes'' corrupted data points during the training.\\

\begin{table}[ht]
\centering
\begin{tabular}{c||r|r||r|r||r|r}
    &  \multicolumn{2}{c||}{1. iteration}  & \multicolumn{2}{c||}{2. iteration}  &  \multicolumn{2}{c}{49. iteration}\\
\hline
                        & corrupted & clean & corrupted & clean & corrupted & clean\\
$u_i$-values            & data points    & data points  & data points    & data points  & data points    & data points\\
\hline\hline
$\gg0$	                &0	        &1  	  &0	      &1	    &1	         &1\\
$\approx 0$  	        &11069      &6251     &565	      &748	    &1905	     &5134\\
$-2.7\cdot 10^{-5}$     &1035       &267      &10504      &5503	    &1*	         &0\\
$-4.0\cdot 10^{-5}$     &0	        &0	      &1035	      &267	    &0	         &0\\
$-5.4\cdot 10^{-5}$     &0	        &0    	  &0	      &0	    &10197	     &1384
\end{tabular}
\caption{Evolution of $u$-vector across 12104 corrupted data points and 6519 clean data points. Note that $1/(12104+6519) = 5.4\cdot 10^{-5}$. Asterisk indicates that this particular $u_i = -1.3\cdot 10^{-5}$.}
\label{tab:u}
\end{table}

\begin{table}[ht]
\centering
\begin{tabular}{l||r|r|r|r|r|r}
    &  \multicolumn{5}{c}{Percentage corrupted training data}\\
\hline
Method                                   & 65\%      & 60\%      & 50\%      & 25\%       & 10\%\\
\hline\hline
Empirical risk minimization              & 0.356-0.389&0.543-0.574&0.666-0.718&0.760-0.780 &0.805-0.807\\
Rockafellian ($\mu = 0.5, \theta = 0.4$)& 0.442   	 &0.610	     &0.694      &0.732       &0.747\\
Rockafellian ($\mu,\theta$ tuned)       & 0.442   	 &0.655	     &0.719      &0.775       &0.785
\end{tabular}
\caption{Test accuracy for empirical risk minimization and Rockafellian relaxation for CIFAR under different levels of corruption.}
\label{tab:CIFARaccuracy}
\end{table}

\noindent {\bf CIFAR Experiments.} We repeat the calculations by constructing a data set from CIFAR by considering three classes (airplane, bird, car) and also leverage Tensorflow 2.10.0. This results in 15000 training data points. We adopt the cross-entropy loss, a fully connected neural net with ReLU activation function, three hidden layers consisting of 320, 320, and 200 units, and a 3-unit output layer utilizing the softmax function. The number of trainable weights in the neural net is 1150040. The image input consists of 32-by-32-by-3 pixels. The computing time for the heuristic (with the same hyperparameters as for MNIST) to carry out 500 epochs and 49 iterations of $u$-optimization is similar to those reported for MNIST. Table \ref{tab:CIFARaccuracy} displays the test accuracy for empirical risk minimization and Rockafellian relaxation for different levels of corruption in a manner parallel to Table \ref{tab:MNISTaccuracy}; the test data set is of size 3000 and is free of corruption. Again, Rockafellian relaxation holds an advantage for higher levels of corruptions. The results for Rockafellian relaxation further improves after a rudimentary tuning of the hyperparameters $\mu$ and $\theta$ across the values $\mu \in \{0.2, 0.3, 0.4, 0.5\}$ and $\theta \in \{0.2, 0.5, 0.8\}$; see the last row in Table \ref{tab:CIFARaccuracy}. As on the MNIST data set, Rockafellian relaxation tends to ``downweight'' data points that are corrupted by assigning them negative $u_i$-values but we omit the details.\\

\noindent {\bf IMDB Experiments.} The third dataset,  IMDB, consists of 50000 movie reviews, each classified as either positive or negative in sentiment. The number of positive and negative reviews are balanced. We generate 25000 training data points of which a fraction are corrupted by incorrect labels: a review with a positive sentiment are assigned a negative sentiment, and vice versa. We tokenize the review text using a WordPiece  subword segmentation algorithm \cite{DBLP:journals/corr/WuSCLNMKCGMKSJL16}. We adopt the DistilBERT architecture \cite{Sanh2019DistilBERTAD} with pre-trained weights made available through the Hugging Face model repository (https://huggingface.co/distilbert-base-uncased). Additionally, we adopt a low-rank adaptation (LoRA) for large language models \cite{hu2022lora} to reduce the number of trainable weights from 67584004 to 628994. The loss function is binary cross-entropy. 

Using Pytorch 2.1.0, the computing time to execute 30 iterations of the heuristic (for a total of 300 epochs) is 960 minutes. This includes 29 updates of $u$, each of which amounts to solving a linear program involving 25000 variables. Thus, the task of updating the $u$-vector is insignificant compared to the overall effort. The capacity of the DistilBERT network is significantly larger than those of the networks used for the MNIST and CIFAR data set. This, combined with a larger data sets results in longer computing times. Again, we use $\mu = 0.5$ and $\theta = 0.4$.

\begin{figure}
\centering
\includegraphics[width=0.99\textwidth]{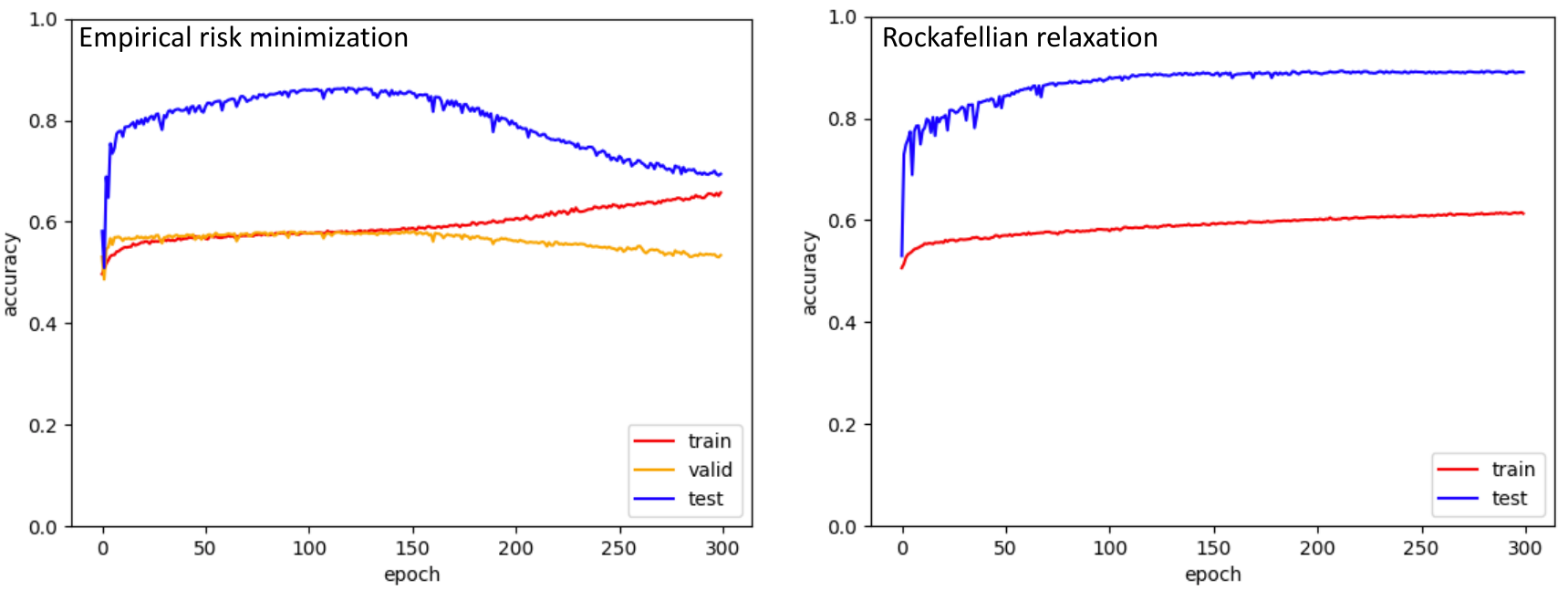}
\caption{Training, validation, and test accuracy for empirical risk minimization (left) and for Rockafellian relaxation (right) on IMDB with 40\% corrupted labels.}\label{IMBD}
\end{figure}

Table \ref{tab:IMBDaccuracy} displays the test accuracy for empirical risk minimization and Rockafellian relaxation at different levels of corruption in a manner parallel to Table \ref{tab:MNISTaccuracy}. The test data set consists of the remaining 25000 reviews and is free of corruption. Again, Rockafellian relaxation holds an advantage for higher levels of corruptions. We notice that the reported accuracy levels for empirical risk minimization rely on an effective means of stopping the algorithm. As seen in Figure \ref{IMBD}(left), the training accuracy (red line) steadily improves under empirical risk minimization while the test accuracy (blue line) starts dropping after 140 epochs. In fact, after 300 epochs the test accuracy has dropped to 0.70. In contrast, the results from Rockafellian relaxation (see Figure \ref{IMBD}(right)) are more stable with the test accuracy (blue line) holding steady at 0.89 and thus making the choice of stopping criterion less important.\\

\begin{table}[ht]
\centering
\begin{tabular}{l||r|r|r|r|r|r}
    &  \multicolumn{5}{c}{Percentage corrupted training data}\\
\hline
Method                                  & 40\%   & 30\%     & 20\%        & 10\%        & 5\%          & 0\%\\
\hline\hline
Emp. risk    &0.853-0.864 &0.898-0.898 &0.903-0.904  &0.910-0.911 &0.913-0.914 & 0.918-0.921\\
Rockafellian &0.893       &0.906       &0.912        &0.917  	  &0.919       &0.923
\end{tabular}
\caption{Test accuracy for empirical risk minimization and Rockafellian relaxation on IMBD under different levels of corruption.}
\label{tab:IMBDaccuracy}
\end{table}

\noindent {\bf Acknowledgement.} This work is supported in part by AFOSR (Math. Optimization) under 21RT0484.

\bibliographystyle{plain}
\bibliography{refs}

\appendix
\section{Appendix}
\renewcommand{\theequation}{\thesection.\arabic{equation}}

\state Proof of Proposition \ref{lemma:: ExactnessGeom}. Since $\mathcal{V}(0) \leq \mathcal{V}(u) - \langle\bar{y}, u\rangle$ for all $u\in\reals^m$ if and only if $\mathcal{V}(0) = \ninf_u \mathcal{V}(u) - \langle\bar{y}, u \rangle = - \mathcal{V}^*(\bar{y})$, the characterization of exactness follows immediately.\eop

\state Proof of Proposition \ref{thm:minChar}. We observe that 
\begin{equation}\label{eqn:proofminChar}
(u^\star, x^\star) \in \nargmin_{u,x} f(u,x) - \langle \bar y, u\rangle \iff u^\star \in \nargmin_u \mathcal{V}(u) - \langle \bar y, u\rangle, \;\; x^\star \in \nargmin_x f(u^\star, x).
\end{equation}
Since $-\cV^*(\bar y) = \inf_u \cV(u) - \langle \bar y, u\rangle$, this implies that $(0, x^\star) \in \nargmin_{u,x} f(u,x) - \langle \bar y, u\rangle$ if and only if $\mathcal{V}(0) = - \mathcal{V}^*(\bar{y})$ and $x^\star \in \nargmin_x \phi(x)$, which in light of Proposition \ref{lemma:: ExactnessGeom} establishes the equivalence of exactness with (\ref{eqn:: ExactnessImplication}).

Next, we establish the necessity of (\ref{eqn:: StrictExactnessImplication}) for strict exactness. If $f$ is strictly exact and $(u^\star, x^\star) \in \nargmin_{u,x} f(u,x) - \langle \bar y, u\rangle$, then by \eqref{eqn:proofminChar} one has $u^\star \in \nargmin_u \mathcal{V}(u) - \langle \bar y, u\rangle$ and $x^\star \in \nargmin_x f(u^\star, x)$. Thus, $\cV(u^\star) - \langle \bar y, u^\star\rangle = \ninf_{u} \cV(u) - \langle \bar y, u\rangle = - \cV^*(\bar y) = \cV(0)$ by Proposition \ref{lemma:: ExactnessGeom}. Definition \ref{def:exactness} then ensures that $u^\star = 0$ and $x^\star \in \argmin_x \phi(x)$.

For the sufficiency, we note that \eqref{eqn:: StrictExactnessImplication} implies via \eqref{eqn:proofminChar} that $\cV(u) \geq \cV(0) + \langle \bar y, u\rangle$ for all $u\in\reals^m$ so $f$ is exact and $- \cV^*(\bar y) = \cV(0)$ by Proposition \ref{lemma:: ExactnessGeom}. Let $\tilde u \in \dom \cV$ be nonzero, then \eqref{eqn:: StrictExactnessImplication} imposes that $(\tilde u, \tilde x) \notin \nargmin_{u,x} f(u,x) - \langle \bar y, u\rangle$ regardless of $\tilde x$. Hence, $f(\tilde u,\tilde x) - \langle \bar{y}, \tilde u \rangle > \ninf_u \mathcal{V}(u) - \langle\bar{y}, u \rangle = - \mathcal{V}^*(\bar{y}) = \cV(0)$. So when $\tilde x\in \argmin_x f(\tilde u,x)$ it is revealed that $\mathcal{V}(\tilde u) > \cV(0) + \langle \bar{y}, \tilde u \rangle$. The same strict inequality holds trivially if $\tilde u\not\in \dom \cV$.\eop

\state Proof of Corollary \ref{cor:SuffCondExactConvex}. This holds by the subgradient inequality of convex analysis.\eop

\state Proof of Proposition \ref{prop:SuffCondExactCalm}. This fact is related to \cite[Theorem 11.61]{VaAn} and follows directly from Definition \ref{def:exactness} because, for $u\in\reals^m$ and $\theta \geq \bar \theta$, one has $\ninf_x f_\theta(u,x) \geq (\theta - \bar\theta)\|u-\bar u\| + \ninf_x f_\theta(\bar u,x)$.\eop

\state Proof of Theorem \ref{thm:asymptExact}. Suppose that the functions $\{f^\nu, \nu\in\nats\}$ are limit-exact Rockafellians supported by $\bar y$. Then, there is an exact Rockafellian $f:\reals^m\times\reals^n\to \Reals$ for the actual problem, supported by $\bar y$, to which the functions $f^\nu$ epi-converge. Consider the functions given by
$\tilde f^\nu(u,x) = f^\nu(u,x) - \langle y^\nu, u\rangle$  and $\tilde f(u,x)  = f(u,x) - \langle \bar y, u\rangle.$
It follows directly from \eqref{eqn:liminfcondition} and \eqref{eqn:limsupcondition} that the epi-convergence of $f^\nu$ to $f$ guarantees that $\tilde f^\nu$ epi-converges to $\tilde f$.

We consider three cases. (i) If $\tilde f$ is proper, then \cite[Theorem 5.5(b)]{primer} applies and we conclude that
\begin{equation}\label{eqn:proofConv}
(\hat u, \hat x) \in \nargmin_{u,x} \tilde f(u,x).
\end{equation}
(ii) If $\tilde f$ takes the value $-\infty$ at one or more points, then the arguments in the proof of \cite[Theorem 5.5(a,b)]{primer} carry over and we again conclude that \eqref{eqn:proofConv} holds. (iii) The case $\tilde f$ being identical to $\infty$ is ruled out by the assumption $\dom \phi \neq \emptyset$. Regardless, \eqref{eqn:proofConv} implies that the necessary optimality condition \eqref{eqn:optcondbasic} holds.

If $\{f^\nu, \nu\in\nats\}$ are strictly limit-exact Rockafellians, then $f$ is strictly exact. Thus, \eqref{eqn:proofConv} and Proposition \ref{thm:minChar} imply that $\hat u = 0$ and $\hat x\in \nargmin_x \phi(x)$.\eop

\end{document}